\documentclass[a4,leqno]{article}
\usepackage{amssymb}
\usepackage{amsfonts}
\usepackage{amsmath}
\usepackage{amsbsy}
\usepackage{amsgen}
\usepackage{amsopn}
\usepackage{amscd}
\usepackage{amstext}
\usepackage{amsxtra}
\usepackage{makeidx}

\title{
Weak conformality of stable stationary maps \\ 
for a functional related to conformality
}
\author{Shigeo Kawai \& Nobumitsu Nakauchi}
\date{}

\begin{document}

\maketitle

\vskip  2ex

\noindent
ABSTRACT. \ \ 
Let $(M,\,g)$, $(N,\,h)$ be compact Riemannian manifolds without boundary, 
and let $f$ be a smooth map from $M$ into $N$. 
We consider a covariant symmetric tensor $T_f$ $=$ 
${\displaystyle f^*h \,-\,\frac{1}{m} \|df\|^2 g}$, 
where $f^*h$ denotes the pull-back metric of $h$ by $f$. 
The tensor $T_f$ vanishes if and only if 
the map $f$ is weakly conformal. 
The norm $\|T_f\|$ is a quantity which is a measure of conformality 
of $f$ at each point. 
We are concerned with 
maps which are critical points of 
the functional 
$\Phi (f)$ $=$ ${\displaystyle \int_M \|T_f\|^2dv_g}$. 
We call such maps {\it C-stationary maps}. 
Any conformal map or more generally 
any weakly conformal map is a C-stationary map. 
It is of interest to find 
when a C-stationary map is a (weakly) conformal map. 

In this paper 
we prove the following result. 
If $f$ is a stable C-stationary maps 
from the standard sphere ${\mathrm S}^m$ $(m \geq 5)$ 
or into the standard sphere ${\mathrm S}^n$ $(n \geq 5)$, 
then $f$ is a weakly conformal map.

\vskip 1ex


\section{Introduction}

Let $(M,\,g)$, $(N,\,h)$ be 
Riemannian manifolds. 
A smooth map $f$ from $M$ into $N$ is a {\it comformal map} 
if and only if there exists a smooth positive function $\varphi$ on $M$ 
such that $f^*h = \varphi g$, 
where $f^*h$ denotes the pullback metric 
of $h$ by $f$, i.e., 
\begin{eqnarray*}
(f^*h)(X,\,Y)\ = \ h\bigl(df(X),\,df(Y)\bigr).
\end{eqnarray*} 
In this situation we utilize
a covariant tensor 
\begin{eqnarray*}
T_f :\ = \ f^*h\ - \ \frac{1}{m}\|df\|^2g.
\end{eqnarray*} 
where
\begin{eqnarray*}
\|df\|^2 \ = \ \sum_{i }
h\bigl(df(e_i),\,df(e_i)\bigr).
\end{eqnarray*} 
Then {\it $f$ is a conformal map 
if and only if $T_f = 0$, 
unless $df \neq 0$}.  
We consider a functional 
\begin{eqnarray*}
\Phi (f) 
= 
\int_M 
\left\|T_f\right\|^2 dv_g, 
\end{eqnarray*}
where $dv_g$ denotes the volume form of $(M,\,g)$, 
and 
\begin{eqnarray*}
\|T_f\|^2 = \sum_{i,\,j} T_f(e_i,\,e_j)^2\,.
\end{eqnarray*}
($e_i$ is a local orthonormal frame on $(M,\,g)$. )
Minimizers of $\Phi$ are close to 
conformal maps, 
even if there does not exist any conformal map 
from $M$ into $N$.  
In \cite{N}, the second author introduced 
the above functional $\Phi$ and 
proved the first variation formula, the second variation formula, 
a quasi-monotonicity formula and a Bochner type formula.
We call a map $f$ {\bf C-stationary} 
if it is a critical point of the functional $\Phi$, 
i.e., if the first variation of $\Phi$ at $f$ vanishes. 
Any conformal map or more generally 
any weakly conformal map is a C-stationary map. 
It is of interest to find 
when a C-stationary map is a (weakly) conformal map. 
In this paper we prove the following two theorems 
for stable C-stationary maps. 

\vskip 3ex 

\noindent 
{\bf Theorem 1}.\ \  
Let $f$ be a stable C-statinary map
from the standard sphere ${\mathrm S}^m$ 
into a Riemannian manifold $N$. 
If $m$ $\geq$ $5$, 
then $f$ is a weakly conformal map.  

\vskip 3ex 

\noindent 
{\bf Theorem 2}.\ \  
Let $f$ be a stable C-statinary map 
from a Riemannian manifold $M$ 
into the standard sphere ${\mathrm S}^n$. 
If $n$ $\geq$ $5$, 
then $f$ is a weakly conformal map.

\vskip 3ex

\noindent 
The contents of this paper 
are as follows: 

\vskip 3ex 

1.  Introduction 

\vskip 1ex 

2.  Weakly conformal maps and the functional $\Phi$  

\vskip 1ex 

3.  Preliminaries

\vskip 1ex 

4.  Stable C-stationary maps from spheres 

\vskip 1ex 

5.  Stable C-stationary maps into spheres

\vskip 5ex

\section{Weakly conformal maps and the functional $\Phi$} 


Let $f$ be a smooth map from a Riemannian manifold $(M,\,g)$ 
into a Riemannian manifold $(N,\,h)$. 
In this section we give a tensor $T_f$ of conformality 
for any smooth map $f$.  
We recall here the following two notions.

\vskip 3ex

\noindent 
{\bf Definition 1}. \ 
(i)\ A map $f$ is {\bf conformal} 
if there exists a smooth {\bf positive} function $\varphi$ on $M$ 
such that 
\begin{eqnarray}
f^*h=\varphi g\,.
\label{eq:conformal}
\end{eqnarray}
(ii)\ A map $f$ is {\bf weakly conformal} 
if there exists a smooth {\bf non-negative} function $\varphi$ on $M$ 
satisfying $(\ref{eq:conformal})$. 

\vskip 3ex

The condition $(\ref{eq:conformal})$ is 
equivalent to 
\begin{eqnarray}
f^*h=\frac{1}{m} \|df\|^2 g\,,
\label{eq:conformal_2}
\end{eqnarray}
since taking the trace of the both sides of $(1)$ 
(w.r.t. the metric $g$), we have $\|df\|^2$ $=$ 
$m \varphi$, i.e., $\varphi$ $=$ $\frac{1}{m}\|df\|^2$. 
Then $f$ is weakly conformal 
if and only if it satisfies $(\ref{eq:conformal_2})$. 
Note that $f$ is weakly conformal if and only if 
for any point $x$ $\in M$, 
$f$ is conformal at $x$ or $df_x = 0$.


Taking the above situation into consideration, 
we utilize the covariant tensor 
\begin{eqnarray*}
T_{f}\ \stackrel{def}{=} f^*h - \frac{1}{m} \|df\|^2 g, 
\end{eqnarray*}
i.e., 
\begin{eqnarray*}
T_{f}(X,\,Y) 
& \stackrel{def}{=} & 
(f^*h)(X,\,Y) - \frac{1}{m} \|df\|^2 g(X,\,Y) 
\\ 
& = & 
h\bigl(df(X),\,df(Y)\bigr) - \frac{1}{m} \|df\|^2 g(X,\,Y)\,,
\end{eqnarray*}
where $f^*h$ denotes the pull-back of the metric $h$.

\vskip 3ex 

\noindent 
{\bf Remark 1}. 
In the case of $m = 2$, 
the quantity $T_f$ is equal to 
the stress energy tensor 
(up to the sign) 
\begin{eqnarray*}
S_f = f^*h - \frac{1}{2} \|df\|^2 g
\end{eqnarray*}
in the harmonic map theory. 
(See Eells-Lemaire \cite{E-L2}, p.392. )

\vskip 3ex 

\noindent 
{\bf Lemma 1}. \ 
\begin{enumerate}
\renewcommand{\labelenumi}{(\alph{enumi})}
\setlength{\itemindent}{30pt} 
\setlength{\labelsep}{10pt}    
\setlength{\itemsep}{0pt}    
\setlength{\parsep}{0pt}    
  \item\ $T_f$ is symmetric, i.e., $T_f(X,\,Y) = T_f(Y,\,X)$. 
  \item\ $f$ is weakly conformal if and only if $T_f = 0$. 
  \item\ $T_f$ is trace-free (with respect to the metric $g$), i.e., 
         \begin{eqnarray*} 
          (g,\,T_f)
          \ = \ 
          \text{Trace}_g T_f
          \ = \ 
          \sum_{i}
                T_f(e_i,\,e_i)
          \ = \ 0\,.
         \end{eqnarray*} 
  \item\ 
          The pairing of 
          the pull-back metric $f^*h$ 
          and the tensor $T_f$ 
         is equal to the norm $\|T_f\|$, i.e.,  
         \begin{eqnarray*} 
          (f^*h,\,T_f)
           \ = \ 
          \sum_{i,\,j}
                (f^*h)(e_i,\,e_j) T_f(e_i,\,e_j)
           \ = \ \|T_f\|^2\,.
         \end{eqnarray*} 
  \item\ ${\displaystyle \|T_f\|^2 = \|f^*h\|^2 - \frac{1}{m} \|df\|^4}$. 
\end{enumerate}
In the above equalities, 
The product $(\ ,\,\ )$ denotes the pairing of the covariant 2-tensors, i.e., 
\begin{eqnarray*}
(A,\,B)\ = \ \sum_{i,\,j\ =\ 1}^n A(e_i,\,e_j) B(e_i,\,e_j) 
\end{eqnarray*}
for any covariant 2-tensors 
$A$, $B$, 
where $e_i$ $(i = 1,\,\cdots\,,\,m)$ is an orthonormal frame.

\vskip 5ex 

\noindent
{\bf Proof}.\ \ 
(a) follows directly from the definition of $T_f$. 

\noindent 
(b) follows easily from the argument mentioned above, i.e., 
\\ \\ 
\qquad\qquad 
$f$ is a weakly conformal 
\\ 
\qquad
$\Leftrightarrow$
\\ 
\qquad\qquad 
There exists a smooth {\it non-negative} function $\varphi$ s.t. 
$f^*h = \varphi g$
\\ 
\qquad
$\Leftrightarrow$
\\ 
\qquad\qquad 
$f^*h = \frac{1}{m}\,\|df\|^2\,g$
\\ 
\qquad
$\Leftrightarrow$
\\ 
\qquad\qquad 
$T_f = 0$\,. 
\vskip 3ex 
\noindent 
\noindent
Here if $f^*h = \varphi g$, 
then taking the trace of the both sides, 
we have $\|df\|^2$ $=$ $m \varphi$, i.e., 
$\varphi$ $=$ ${\displaystyle \frac{1}{m}}$ $\|df\|^2$.   

\vskip 3ex 

\noindent 
(c)\ \  
Note $(g,\,T_f)$ $=$ $\text{Trace}_g T_f$. 
Moreover we have 
\begin{eqnarray*}
\text{Trace}_g T_f
& = & 
\sum_{i}
   T_f(e_i,\,e_i)
\\ 
& = & 
\sum_{i}
    \left\{h\bigl(df(e_i),\,df(e_i)\bigr) - \frac{1}{m} \|df\|^2 g(e_i,\,e_i) \right\}
\\ 
& = & 
\sum_{i}
    h\bigl(df(e_i),\,df(e_i)\bigr)
\ - \ 
\frac{1}{m} \|df\|^2
\sum_{i}
    g(e_i,\,e_i)
\\ 
& = & 
\sum_{i} h\bigl( df(e_i),\,df(e_i)\bigr)
\ - \ 
\|df\|^2
\\
& = & 
\|df\|^2
\ - \ 
\|df\|^2
\\ 
& = & 
0\,. 
\end{eqnarray*}


\begin{align}
(f^*h,\,T_f)
& \ = \ 
\bigl(T_f + \frac{1}{m} \|df\|^2g,\,T_f\bigr)
\tag{d} \\ 
& \ = \  
\bigl(T_f,\,T_f\bigr) 
\ - \ 
\frac{1}{m} \|df\|^2
\bigl(g,\,T_f\bigr)
\notag \\ 
& \ = \  
\|T_f\|^2\,. 
\notag 
\end{align}


\begin{align}
\|T_f\|^2\ 
& \ = \ 
(f^*h,\,T_f)  
\tag{e}
\\ 
& \ = \ 
\bigl(f^*h,\,f^*h - \frac{1}{m} \|df\|^2g\bigr)
\notag \\ 
& \ = \ 
\|f^*h\|^2
\ - \ 
\frac{1}{m} \|df\|^2 
\bigl(f^*h,\,g\bigr)
\notag \\ 
& \ = \ 
\|f^*h\|^2
\ - \ 
\frac{1}{m} \|df\|^4\,. 
\notag 
\end{align}
Thus we obtain Lemma 1. $\square$

\vskip 3ex

In this paper, 
we are concerned with the functional of the norm of $T_f$ 
\begin{eqnarray*}
\Phi (f)\ =\ {\displaystyle \int_M \|T_f\|^2dv_g}\,.
\end{eqnarray*} 
This quantity $\Phi (f)$ gives 
a measure of the conformality of 
maps $f$. 
Note that if $f$ is a conformal map, 
then $\Phi (f)$ vanishes.

\vskip 3ex

\section{Preliminaries}

In this section
we give a technical lemma (Lemma 2), 
the first variation formula and the second variation formula 
with some notations and definitions. 
The two formulas are obtained in the second author's paper \cite{N-T}. 
For reader's convenience we give their proofs of these two formula 
in the appendix at the end of the paper.

Take any smooth deformation $F$ of $f$, i.e., 
any smooth map 
\begin{eqnarray*}
F\ :\ (-\varepsilon,\,\varepsilon) \times M\ \longrightarrow\ N
\ \ \text{s.t.}\ \ F(0,\,x) = f(x), 
\end{eqnarray*}
where $\varepsilon$ is a positive constant. 
Let $f_t(x) = F(t,\,x)$. Then we have $f_0(x) = f(x)$. 
We often say a deformation $f_t(x)$ 
instead of a deformation $F(t,\,x)$. 
Let $X$ $=$ $\left.dF({\textstyle \frac{\partial}{\partial t}})
\right|_{t=0}$ denotes 
the variation vector field of the deformation $F$. 

We define an $f^{-1} TN$-valued 1-form $\sigma_f$ on $M$ by
\begin{eqnarray}
\sigma_f (X) = \sum_j T_f (X,\,e_j) df(e_j)
\label{eq:definition_of_sigma}
\end{eqnarray}
where $\{e_j\}$ is an orthonormal frame.
The 1-form $\sigma_f$ plays an important role in our arguments.

\vskip 3ex



We first give the following lemma, 
which we often use in our arguments.

\vskip 3ex

\noindent
{\bf Lemma 2}.\ \ 
\begin{eqnarray}
\sum_j 
  h\bigl(Z,\,dF(e_j)\bigr)\,
  T_F\bigl(W,\,e_j\bigr)
\ = \ 
  h\bigl(Z,\,\sigma_F (W)\bigr)\,.
\label{eq:1_in_Lemma2}
\end{eqnarray}
In particular 
\begin{eqnarray}
\sum_j 
  h\bigl(Z,\,df(e_j)\bigr)\,
  T_f\bigl(W,\,e_j\bigr)
\ = \ 
  h\bigl(Z,\,\sigma_f (W)\bigr)\,.
\label{eq:2_in_Lemma2}
\end{eqnarray}
\vskip -3ex 
\begin{eqnarray}
\|T_f\|^2 
\ = \ 
 \sum_i h\bigl(df(e_i),\,\sigma_f (e_i)\bigr)\,.
\label{eq:3_in_Lemma2}
\end{eqnarray}

\vskip 3ex

\noindent
{\bf Proof of Lemma 2}. 
The equality $(\ref{eq:1_in_Lemma2})$ 
easily follows from the definition of $\sigma_F$. 
Indeed, since $h(A,\,B)T_F(C,\,D)$ $=$ $h(A,\,T_F(C,\,D)B)$, 
we have
\begin{eqnarray*}
\sum_j 
  h\bigl(Z,\,dF(e_j)\bigr)\,
  T_F\bigl(W,\,e_j\bigr)
& = & 
  h\bigl(Z,\,\sum_j T_F\bigl(W,\,e_j) dF(e_j)\bigr)\,
\\ 
& = & 
  h\bigl(Z,\,\sigma_F (W)\bigr)\,.
\end{eqnarray*}
The equality $(\ref{eq:2_in_Lemma2})$ follows from 
$(\ref{eq:1_in_Lemma2})$ at $t =0$. 
Furthermore let $Z$ $=$ $df(e_i)$ and let $W$ $=$ $e_i$ 
in $(\ref{eq:2_in_Lemma2})$, 
and sum with respect to $i$. 
Then we have 
\begin{eqnarray*}
\sum_{i,\,j} h\bigl(df(e_i),\,df(e_j)\bigr)\,T_f(e_i,\,e_j)
\ = \ 
 \sum_i h\bigl(df(e_i),\,\sigma_f (e_i)\bigr)\,.
\end{eqnarray*}
Since $(f^*h)(e_i,\,e_j)$ $=$ $h\bigl(df(e_i),\,df (e_j)\bigr)$,  
the above equality and Lemma 1 (d) imply 
$(\ref{eq:3_in_Lemma2})$. 
$\square$

\vskip 3ex

The first variation of the functional $\Phi$ is given 
by the following formula which is given and proved in \cite{N}.

\vskip 3ex 

\noindent
{\bf Proposition 1} (first variation formula).\ \ 
\begin{eqnarray*}
\left.\frac{d\Phi(f_t)}{dt}  \right|_{t=0}
& = & 
- 4
 \int_M 
h\left( 
X
,\,
\mbox{\text div}_g \sigma_f
\right) \,dv_g\,,
\end{eqnarray*}
where 
$\mathrm{div}_g\sigma_f$ denotes the divergence of $\sigma_f$, 
i.e., $\mathrm{div}_g\sigma_f$ $=$ 
${\displaystyle \sum_{i=1}^m (\nabla_{e_i}\sigma_f)(e_i)}$. 

%

\vskip 3ex

We give here the notion of C-stationary maps.  


\vskip 3ex

\noindent
{\bf Definition 2}.\ \ 
We call a smooth map $f$ 
{\bf C-stationary} 
if the first variation of $\Phi ( f )$ (at $f$)
identically vanishes, i.e., 
\begin{eqnarray*}
\left.\frac{d\Phi(f_t)}{dt}  \right|_{t=0}
& = & 
0
\end{eqnarray*}
for any smooth deformation $f_t$ of $f$. 
By Proposition 1, 
a smooth map $f$ is {\it C-stationary} 
if and only if it satisfies the equation 
\begin{eqnarray}
\mathrm{div}_g \sigma_f \ = \ 0\,,
\label{equation:Euler-Lagrange equation}
\end{eqnarray}
which is called the Euler-Lagrange equation 
for the functional $\Phi ( f )$, 
where $\sigma_f$ is the covariant tensor 
defined by $(\ref{eq:definition_of_sigma})$.

\vskip 3ex

We give the second variation formula 
for the functional $\Phi ( f )$. 
Take any smooth deformation $F$ of $f$ 
with two parameters, i.e., 
any smooth map 
\begin{eqnarray*}
F\ :\ (-\varepsilon,\,\varepsilon) \times (-\delta,\,\delta) \times M
\ \longrightarrow\ N
\ \ \text{s.t.}\ \ F(0,\,0,\,x) = f(x)\,. 
\end{eqnarray*}
Let $f_{s,\,t}(x) = F(s,\,t,\,x)$, 
and we often say a deformation $f_{s,\,t}(x)$ 
instead of a deformation $F(s,\,t,\,x)$. 
Let 
\begin{center}
$X$ $=$ $\left.dF(\frac{\partial}{\partial s})\right|_{s,\,t=0}$,\ \
$Y$ $=$ $\left.dF(\frac{\partial}{\partial t})\right|_{s,\,t=0}$\,
\end{center}
denote 
the variation vector fields of the deformation $f_{s,\,t}$. 
Then we have the following second variation formula. 

\vskip 3ex 

\noindent
{\bf Proposition 2} (second variation formula).\ \ 
\begin{eqnarray}
\lefteqn{
\frac{1}{4}
\left.\frac{\partial^2\Phi(f_{s,\,t})}{\partial s \partial t}\right|_{s,\,t = 0}
\ = \ 
\int_M  
h\bigl(
\mbox{\text Hess}_{F}
\bigl(
{\textstyle \frac{\partial}{\partial s}}
,\,
{\textstyle \frac{\partial}{\partial t}}
\bigr),\,
\mathrm{div}_g\sigma_{f}
\bigr)
dv_g
} 
\label{eq:second variation formula}
\\ 
& &
\ + \ 
\int_M
\sum_{i,\,j} 
h\bigl(
\nabla_{e_i} X,\,
\nabla_{e_j} Y \bigr)\,
T_f(e_i,\,e_j)\,
dv_g
\nonumber \\ 
&  & 
\ + \ 
\int_M
\sum_{i,\,j} 
h\bigl(\nabla_{e_i} X,\,df(e_j) \bigr)\,
h\bigl(\nabla_{e_i} Y,\,df(e_j)\bigr) 
dv_g
\nonumber \\ 
&  & 
\ + \ 
\int_M
\sum_{i,\,j} 
h\left(\nabla_{e_i} X,\,df(e_j)\right)\, 
h\left(df(e_i),\,\nabla_{e_j} Y\right) 
dv_g
\nonumber \\ 
&  & 
\ - \ 
\frac{2}{m}
\int_M
\sum_{i} 
h\left(\nabla_{e_i} X,\,df(e_i)\right)\, 
\sum_{j} 
h\left(\nabla_{e_j} Y,\,df(e_j)\right)\, 
dv_g
\nonumber \\ 
&  & 
\ + \ 
\int_M
\sum_{i,\,j} 
h\bigl({}^N\!
R\left(
X,\,
df(e_i)\,
\right)Y,\,
df(e_j)
\bigr)\,
T_f(e_i,\,e_j)\,
dv_g\,,
\nonumber 
\end{eqnarray}
where 
$\mathrm{Hess}_f$ denotes the Hessian of $f$, i.e., 
$\mathrm{Hess}_f (Z,\,W)$ $=$ $(\nabla_Z df)(W)$ $=$ $(\nabla_W df)(Z)$\,.

\vskip 3ex

\noindent
{\bf Remark 2}.\ 
Note that the first term in the right hand side 
vanishes if $f$ is a C-stationary map. 

\vskip 3ex

\noindent
{\bf Remark 3}.\ 
The last term of the right hand side 
in Proposition 2 is equal to 
\begin{eqnarray*}
\int_M
\sum_{i} 
h\bigl({}^N\!
R\left(
X,\,
df(e_i)\,
\right)Y,\,
\sigma_f (e_i)
\bigr)\,
dv_g\,,
\end{eqnarray*}
since 
\begin{eqnarray*}
\lefteqn{
\sum_{i,\,j} 
h\bigl({}^N\!
R\left(
X,\,
df(e_i)\,
\right)Y,\,
df(e_j)
\bigr)
T_f(e_i,\,e_j)\,
} \nonumber \\ 
& = & 
\sum_{i} 
h\bigl({}^N\!
R\left(
X,\,
df(e_i)\,
\right)Y,\,
\sum_j 
T_f(e_i,\,e_j)\,
df(e_j)
\bigr)\,.
\nonumber \\ 
& = & 
\sum_{i} 
h\bigl({}^N\!
R\left(
X,\,
df(e_i)\,
\right)Y,\,
\sigma_f
(e_i)
\bigr)\,.
\nonumber
\end{eqnarray*}

\vskip 3ex

\vskip 3ex

\noindent
{\bf Definition 3}.\ \ 
We call a C-stationary map {\it stable} 
if the second variation of $\Phi ( f )$ (at $f$)
is nonnegative, i.e., 
\begin{eqnarray*}
\left.\frac{d^2\Phi(f_t)}{dt^2}  \right|_{t=0}
& \geq & 
0.
\end{eqnarray*}

\vskip 3ex

We give an example of the stable C-stationary map 
which is not weakly conformal.

\vskip 3ex 

\noindent
{\bf Example}. \ \  
Let us define a map
$$
	f: M = S^1 \times S^1 \times \cdots \times S^1 \to 
		N = S^1_k \times S^1 \times \cdots \times S^1
$$
by
$$
	f(x_1, x_2 \ldots x_\ell) = (k x_1, x_2 \ldots x_\ell)
$$
where $S^1$ (resp. $S^1_k$) denotes the sphere of dimension 1 with
radius 1 (resp. $k$) centered at the origin of $\mathbb R^2$. 
Obviously $f$ is not weakly conformal if $k \neq 1$. 
Let us take orthonormal unit parallel vector fields 
$e_1, e_2 \ldots e_\ell$ on $M$ where $e_i$ is tangent to the
$i$-th factor, and $E_1, E_2 \ldots E_\ell$ on   $N$ where $E_j$
is tangent to the $j$-th factor. Hence we have
$$
	df (e_1) = k E_1, \quad df (e_i) = E_i \quad (i \neq 1).
$$
Every vector field $X$ on $N$ is written as $X = \sum_i \varphi_i E_i$
where $\varphi_i$ is a function on $N$. 

It is easy to show
$$
	T_f (v,w) = C \langle v, w \rangle
$$
where 
$$
	C = \begin{cases}
		\frac{(k^2 - 1)(\ell - 1)}{\ell} & \mbox{
			($v$ and $w$ are tangent to the first factor)} \cr
		- \frac{k^2 - 1}{\ell} & \mbox{
			(otherwise)} \cr
            \end{cases}
$$
Hence we obtain for a tangent vector
$v$ on $M$,          
$$
	\sigma_f (v) = \frac{k(k^2 - 1)(\ell - 1)}{\ell} \langle v, e_1
			\rangle  E_1
			- \frac{k^2 - 1}{\ell} \sum_{j = 2}^\ell 
				\langle v, e_j \rangle E_j,
$$ 
which implies that $ {\rm div\,} \sigma_f = 0$, i.e., the map $f$ is C-stationary.

We have only to prove that $f$ is stable, i.e., 
the second variation at $f$ is non-negative. 
Denoting $\varphi_i \circ f$ by $\psi_i$, 
we culculate terms of
the right hand side of the second variation formula as follows:  
\begin{eqnarray*}
&&	\int_M \sum_{i,j} \langle\nabla_{e_i}X ,  \nabla_{e_j} X \rangle T_f(e_i,e_j) \\
&&	= \int_M \frac{(k^2 - 1)(\ell - 1)}{\ell} 
		\langle \nabla_{e_1} X, \nabla_{e_1} X \rangle
		- \frac{k^2 - 1}{\ell} \sum_{i \geq 2} \langle\nabla_{e_i}X ,
		\nabla_{e_i}X  \rangle \\
&&	= \int_M \frac{(k^2 - 1)(\ell - 1)}{\ell} \sum_j (\nabla_{e_1} \psi_j)^2
		- \frac{k^2 - 1}{\ell} \sum_{i \geq 2} \sum_j (\nabla_{e_i}
		\psi_j)^2 ,
\end{eqnarray*}
$$
	 \int_M \sum_{i,j} \langle\nabla_{e_i}X , df (e_j) \rangle
	\langle \nabla_{e_i} X, df (e_j) \rangle  
	= \int_M \sum_i k^2 (\nabla_{e_i} \psi_1)^2 
		+ \sum_i \sum_{j \geq 2}(\nabla_{e_i} \psi_j)^2,
$$
\begin{eqnarray*}
	 \int_M \sum_{i,j} \langle\nabla_{e_i}X , df (e_j) \rangle
	  \langle df (e_i), \nabla_{e_j} X     \rangle  
	&=& \int_M k^2 (\nabla_{e_1} \psi_1)^2 \\
		+ 2 k \sum_{i \geq 2}  (\nabla_{e_1} \psi_i)
			(\nabla_{e_i}\psi_1) 
	&+& \sum_{i \geq 2, j \geq 2}
		(\nabla_{e_i} \psi_j)(\nabla_{e_j}\psi_i),
\end{eqnarray*}
\begin{eqnarray*}
	- \frac{2}{\ell} \int_M \left(\sum_i \langle \nabla_{e_i} X, df (e_i)
		\rangle \right)^2
	&=& - \frac{2}{\ell} \int_M
		k^2 (\nabla_{e_1} \psi_1)^2 \\
	+ 2k \sum_{i \geq 2} (\nabla_{e_1}\psi_1)(\nabla_{e_i} \psi_i)
	&+& \sum_{i \geq 2}\sum_{ j \geq 2} (\nabla_{e_i}\psi_i)(\nabla_{e_j} \psi_j).
\end{eqnarray*}

Note that the following identities hold:
\begin{eqnarray*}
(\nabla_{e_i} \psi_j)(\nabla_{e_j} \psi_i) &=&
	\nabla_{e_i}(\psi_j \nabla_{e_j} \psi_i) - \psi_j (\nabla_{e_i}
	\nabla_{e_j} \psi _i) \\
(\nabla_{e_i} \psi_i)(\nabla_{e_j} \psi_j) &=&
	\nabla_{e_j}(\psi_j \nabla_{e_i} \psi_i) - \psi_j (\nabla_{e_j}
	\nabla_{e_i} \psi _i).
\end{eqnarray*}
Hence exchanging the orders of iterated integrals, we obtain for each $i$ 
and $j$,
\begin{eqnarray*}
\int_{S^1 \times \cdots \times S^1} 
	(\nabla_{e_i} \psi_j)(\nabla_{e_j} \psi_i) &=&
- \int_{S^1 \times \cdots \times S^1}  \psi_j (\nabla_{e_i}
	\nabla_{e_j} \psi _i), \\
\int_{S^1 \times \cdots \times S^1} 
(\nabla_{e_i} \psi_i)(\nabla_{e_j} \psi_j) &=&
- \int_{S^1 \times \cdots \times S^1} \psi_j (\nabla_{e_j}
	\nabla_{e_i} \psi _i).
\end{eqnarray*}
Thus these integrals coincide because $[e_i, e_j]= 0$. Especially we get
$$
	\int_M ( \nabla_{e_1} \psi_1)( \nabla_{e_i} \psi_i)
	= 
	\int_M ( \nabla_{e_1} \psi_i)( \nabla_{e_i} \psi_1).
$$

Let us denote the function $\nabla_{e_i} \psi_j$ by $a_{ij}$ for simplicity.
From the fact above, 
the right hand side of the second variation formula is equal to
\begin{eqnarray}
& &	\int_M \biggl (	\frac{(k^2 - 1)(\ell - 1)}{\ell}\sum_j (\nabla_{e_1}\psi_j)^2
	- \frac{k^2 - 1}{\ell} \sum_{i \geq 2} \sum_j (\nabla_{e_i}\psi_j)^2 
\label{eq:calculation of example} \\
& &	+ k^2 \sum_i (\nabla_{e_i}\psi_1)^2 +\sum_i \sum_{j \geq 2}
	(\nabla_{e_i}\psi_j)^2 
	+ k^2(\nabla_{e_1}\psi_1)^2  
\nonumber \\
& & + 2k \sum_{i \geq 2}
	(\nabla_{e_1}\psi_i)(\nabla_{e_i}\psi_1) 
	+ \sum_{i \geq 2} \sum_{j \geq 2}
	(\nabla_{e_i}\psi_j)(\nabla_{e_j}\psi_i) 
\nonumber \\
& &	-\frac{2}{\ell} k^2(\nabla_{e_1}\psi_1) -  \frac{4k}{\ell} \sum_{i \geq 2}
	(\nabla_{e_1}\psi_1)(\nabla_{e_i}\psi_i) 
	- \frac{2}{\ell}  \sum_{i \geq 2} \sum_{j \geq 2}
	(\nabla_{e_i}\psi_i)(\nabla_{e_j}\psi_j) \biggr ) 
\nonumber \\
& = & \int_M \biggl (\frac{\ell - 1}{\ell} (3k^2 - 1) a_{11}^2 + \sum_{j \geq 2}
	 \frac{1 - k^2}{\ell} \sum_{i \geq2, j \geq 2}  a_{ij}^2 
	+ \sum_{i \geq 2} 2k\left(1 - \frac{2}{\ell} \right) a_{1i} a_{i1}
\nonumber \\
& &	+ \sum_{j \geq 2} \frac{1}{\ell} \left( k^2 (\ell - 1) + 1 \right) a_{1j}^2 
 + \sum_{i \geq 2} \frac{1}{\ell} \left( k^2(\ell - 1) + 1 \right) a_{i1}^2 
\nonumber \\
& & + \biggl \{  \sum_{i \geq 2, j \geq 2}a_{ij}^2 
	+ \sum_{i \geq 2, j \geq 2} a_{ij} a_{ji}
	 - \frac{2}{\ell} \sum_{i \geq 2, j \geq 2} a_{ii} a_{jj} \biggr \} \biggr ).
\nonumber 
\end{eqnarray}
Note that 
\begin{eqnarray*}
	\sum_{i \geq 2} 2k\left(1 - \frac{2}{\ell}\right) a_{1i}a_{i1} 
	&\geq& - k\left(1 - \frac{2}{\ell}\right)\sum_{i \geq 2}
		\left( a_{1i}^2 + a_{i1}^2 \right) \\
	&=& - k\left(1 - \frac{2}{\ell}\right)\sum_{i \geq 2} a_{1i}^2
		- k\left(1 - \frac{2}{\ell}\right)\sum_{i \geq 2} a_{i1}^2.
\end{eqnarray*}
If $k$ is close to 1, then we have
$$
	 \frac{1}{\ell} \left( k^2 (\ell - 1) + 1 \right) 
	> k\left( 1 - \frac{2}{\ell}\right).
$$
Thus the sum of the third, fourth and fifth terms of the equation
(*) is non-negative.

Since $ \|A \|^2 + {\rm Tr}(A^2) - \frac{2}{\ell - 1}({\rm Tr} A)^2 \geq 0$
for every $(\ell - 1) \times (\ell - 1)$ matrix $A$, we have
$$
   \sum_{i \geq 2, j \geq 2}a_{ij}^2 
	+ \sum_{i \geq 2, j \geq 2} a_{ij} a_{ji}
	 - \frac{2}{\ell} \sum_{i \geq 2, j \geq 2} a_{ii} a_{jj} \geq 0.
$$	
Hence if $k \leq 1$ and $k$ is close to 1, then
the second variation at $f$ is non-negative.

\vskip 3ex 

As is seen, there exist stable C-stationary maps which are
not weakly conformal. We shall show in the next section that this is
not the case when the domain or the range is a standard sphere.

\vskip 3ex

\section{Stable C-stationary maps from spheres}

In this section, we prove Theorem 1.

\vskip 3ex

\noindent
{\bf Proof of Theorem 1}.\ \
Since the standard sphere $\mathrm{S}^m$ is
a submanifold of the Euclidean space $\mathbb{R}^{m+1}$,
we may consider that
the tangent space of the standard sphere $\mathrm{S}^m$
at $x$ $\in$ $\mathrm{S}^m$
is a subspace of the linear space $\mathbb{R}^{m+1}$ $\simeq$
$\mathrm{T}_x\mathbb{R}^{m+1}$.
Let $p=p_x$ denote the canonical projection
from $\mathbb{R}^{m+1}$
onto the tangent space $\mathrm{T}_x\mathbb{S}^{m}$.
Let $E$ be a unit parallel vector field 
on $\mathbb{R}^{m+1}$, 
and let $Z$ be the vector field which is 
the image by the projection $p$ of $E$, i.e.,
\begin{eqnarray*}
Z_x
\,=\
p_{x}(E) 
\end{eqnarray*}
for $x$ $\in$ $\mathrm{S}^m$.
Then we can verify 
\begin{eqnarray}
\nabla_{e_i}Z
\ = \
-\,\varphi\,e_i\,, 
\label{eq:diff_vectorfield_02} 
\nonumber 
\end{eqnarray}
where 
\begin{eqnarray*}
\varphi 
\ = \
\langle E,\,\nu \rangle 
\nonumber 
\end{eqnarray*}
(the notation $\langle\ \ ,\ \ \rangle$ denotes the inner product on $\mathbb{R}^{m+1}$)
and $\nu$ is the unit outer normal vector field on $\mathrm{S}^m$ 
in $\mathbb{R}^{m+1}$. 
Then we have
\begin{eqnarray}
\nabla_{e_i}\bigl(df(Z)\bigr)
& = & 
\bigl(\nabla_{e_i}df\bigr) (Z) 
\ + \ 
df\bigl(\nabla_{e_i}Z \bigr) 
\label{eq:equality_for_vector_field}
 \\ 
& = & 
\bigl(\nabla_{e_i}df\bigr) (Z) 
\ - \ 
\varphi df(e_i) 
\nonumber 
\end{eqnarray}

Take orthonormal parallel vector fields 
$E_1$,\,$\cdots$\,,\,$E_{m+1}$ on $\mathbb{R}^{m+1}$, 
and set $Z_k$ $=$ $p(E_k)$ $(k =1,\,\cdots\,,\,m+1)$. 
Then by $(\ref{eq:equality_for_vector_field})$, 
we see 
\begin{eqnarray}
\nabla_{e_i}\bigl(df(Z_k)\bigr)
& = & 
\bigl(\nabla_{e_i}df\bigr) (Z_k) 
\ - \ 
\varphi_k df(e_i) 
\label{eq:equality_for_vector_field_k}
\end{eqnarray}
where 
\begin{eqnarray*}
\varphi_k 
\ = \
\langle E_k,\,\nu \rangle \,.
\nonumber 
\end{eqnarray*}

The stability of the C-stationary map $f$ implies
the inequality
\begin{eqnarray*}
{\mathrm L}(df(Z_k),\,df(Z_k))
\ \geq\ 0\,,\ \ \text{hence}\ \ 
\sum_{k=1}^{m+1}
{\mathrm L}(df(Z_k),\,df(Z_k))
\geq\ 0\,,
\end{eqnarray*}
where
$\mathrm{L}(X,\,Y)$ denotes 
the right hand side of
the second variation formula $(\ref{eq:second variation formula})$.
We calculate ${\mathrm L}(df(Z_k),\,df(Z_k))$. 
By the definition of $\mathrm{L}(X,\,Y)$,
we see
\begin{eqnarray}
\quad\quad 
\lefteqn{
{\mathrm L}(df(Z_k),\,df(Z_k))
} \label{eq:stability___01}
\\
& = &
\int_M
\sum_{i,\,j}
h \bigl(
\nabla_{e_i}\bigl( df(Z_k)\bigr )
,\,
\nabla_{e_j}\bigl( df(Z_k)\bigr )
\bigr)
\,
T_f (e_i,\,e_j)
\,dv_g
\nonumber \\
&  &
\ + \
\int_M
\sum_{i,\,j}
h \bigl(
\nabla_{e_i}\bigl( df(Z_k)\bigr )
,\,
df(e_j)
\bigr)
\,
h \bigl(
\nabla_{e_i}\bigl( df(Z_k)\bigr )
,\,
df(e_j)
\bigr)
\,dv_g
\nonumber \\
&  &
\ + \
\int_M
\sum_{i,\,j}
h \bigl(
\nabla_{e_i}\bigl( df(Z_k)\bigr )
,\,
df(e_j)
\bigr)
\,
h \bigl(
df(e_i)
,\,
\nabla_{e_j}\bigl( df(Z_k)\bigr )
\bigr)
\,dv_g
\nonumber \\
&  & 
\ - \ 
\frac{2}{m}
\int_M 
\sum_{i} 
h \bigl( 
\nabla_{e_i}\bigl( df(Z_{k})\bigr )
,\,
df(e_i)
\bigr)
\, 
\sum_{j} 
h \bigl( 
\nabla_{e_j}\bigl( df(Z_{k})\bigr )
,\,
df(e_j)
\bigr)
\,dv_g
\nonumber \\ 
&  &
\ + \
\int_M
\sum_{i,\,j}
h\bigl({}^N\!
R\left(
df(Z_k),\,
df(e_i)
\right)
df(Z_k),\,df(e_j)
\bigr)\,
T_f (e_i,\,e_j\bigr)
dv_g .
\nonumber 
\end{eqnarray}
We denote terms on the right hand side by A, B, C, D and E respectively.
By (\ref{eq:equality_for_vector_field_k}) and Lemma 1 (a), (e), 
we have 
\begin{eqnarray}
\mathrm{A}
& = &
\int_M
\sum_{i,\,j}
h \bigl(
\bigl( \nabla_{e_i}df\bigr )(Z_k)
,\,
\bigl( \nabla_{e_j}df\bigr )(Z_k)
\bigr)
\,
T_f (e_i,\,e_j)
\,dv_g
\label{eq:AAA}
\\
&  &
\ - \
2\int_M
\varphi_k
\sum_{i,\,j}
h \bigl(
(\nabla_{e_i}df)(Z_k)
,\,
df(e_j)
\bigr)
T_f (e_i,\,e_j\bigr)
\,dv_g
\nonumber \\
&  &
\ + \
\int_M
\varphi_k ^2
\|T_f\|^2
\,dv_g\,.
\nonumber
\end{eqnarray}
By (\ref{eq:equality_for_vector_field_k}) and Lemma 1 (a), 
we get
\begin{eqnarray}
\mathrm{B}
& = &
\int_M
\sum_{i,\,j}
h \bigl(
\bigl(\nabla_{e_i} df\bigr )(Z_k)
,\,
df(e_j)
\bigr)
\,
h \bigl(
\bigl(\nabla_{e_i} df\bigr )(Z_k)
,\,
df(e_j)
\bigr)
\,dv_g
\label{eq:BBB} \\
&  &
\ - \
2\int_M
\varphi_k
\sum_{i,\,j}
h \bigl(
(\nabla_{e_i}df)(Z_k)
,\,
df(e_j)
\bigr)
h \bigl(
df(e_i)
,\,
df(e_j)
\bigr)
\,dv_g
\nonumber \\
&  &
\ + \
\int_M
\varphi_k ^2
\|f^*h\|^2
\,dv_g\,,
\nonumber \\ 
\mathrm{C}
& = &
\int_M
\sum_{i,\,j}
h \bigl(
\bigl(\nabla_{e_i}df\bigr )(Z_k)
,\,
df(e_j)
\bigr)
\,
h \bigl(
df(e_i)
,\,
\bigl(\nabla_{e_j}df\bigr )(Z_k)
\bigr)
\,dv_g
\label{eq:CCC} \\
&  &
\ - \
2\int_M
\varphi_k
\sum_{i,\,j}
h \bigl(
(\nabla_{e_i}df)(Z_k)
,\,
df(e_j)
\bigr)
h \bigl(
df(e_i)
,\,
df(e_j)
\bigr)
\,dv_g
\nonumber \\
&  &
\ + \
\int_M
\varphi_k ^2
\|f^*h\|^2
\,dv_g\,.
\nonumber
\end{eqnarray}
and
\begin{eqnarray}
\mathrm{D}
& = & 
-\,\frac{2}{m} 
\int_M 
\sum_{i} 
h \bigl( 
(\nabla_{e_i} df)(Z_k)
,\,
df(e_i)
\bigr)
\, 
\sum_{j} 
h \bigl( 
(\nabla_{e_j} df)(Z_k)
,\,
df(e_j)
\bigr)
\,dv_g
\label{eq:DDD} \\
&  & 
\ + \ 
\frac{4}{m} 
\int_M 
\varphi_k  
\sum_{i} 
h \bigl( 
(\nabla_{e_i} df)(Z_k)
,\,
df(e_i)
\bigr)
\, 
\|df\|^2 
\,dv_g
\nonumber \\ 
&  & 
\ - \ 
\frac{2}{m} 
\int_M 
\varphi_k^2 \, 
\|df\|^4
\,dv_g
\nonumber 
\end{eqnarray}
Then by $(\ref{eq:stability___01})$, $(\ref{eq:AAA})$, $(\ref{eq:BBB})$, $(\ref{eq:CCC})$, $(\ref{eq:DDD})$ and Lemma 1 (c), 
we have
\begin{eqnarray}
\lefteqn{
{\mathrm L}(df(Z_k),\,df(Z_k))
} \label{eq:stability___02}
\\
& = &
\int_M
\sum_{i,\,j}
h \bigl(
\bigl( \nabla_{e_i}df(Z_k)\bigr )
,\,
\bigl( \nabla_{e_j}df(Z_k)\bigr )
\bigr)
\,
T_f (e_i,\,e_j)
\,dv_g
\nonumber \\
&  &
\ + \
\int_M
\sum_{i,\,j}
h \bigl(
\bigl( \nabla_{e_i}df(Z_k)\bigr )
,\,
df(e_j)
\bigr)
\,
h \bigl(
\bigl( \nabla_{e_i}df(Z_k)\bigr )
,\,
df(e_j)
\bigr)
\,dv_g
\nonumber \\
&  &
\ + \
\int_M
\sum_{i,\,j}
h \bigl(
\bigl( \nabla_{e_i}df(Z_k)\bigr )
,\,
df(e_j)
\bigr)
\,
h \bigl(
df(e_i)
,\,
\bigl( \nabla_{e_j}df(Z_k)\bigr )
\bigr)
\,dv_g
\nonumber \\
&  & 
-\,\frac{2}{m} 
\int_M 
\sum_{i} 
h \bigl( 
(\nabla_{e_i} df)(Z_k)
,\,
df(e_i)
\bigr)
\, 
\sum_{j} 
h \bigl( 
(\nabla_{e_j} df)(Z_k)
,\,
df(e_j)
\bigr)
\,dv_g
\nonumber \\ 
&  &
\ + \
\int_M
\sum_{i,\,j}
h\bigl({}^N\!
R\left(
df(Z_k),\,
df(e_i)
\right)
df(Z_k),\,df(e_j)
\bigr)\,
T_f (e_i,\,e_j)
\,dv_g .
\nonumber \\
&  &
\ - \
6\int_M
\varphi_k
\sum_{i,\,j}
h \bigl(
(\nabla_{e_i}df)(Z_k)
,\,
df(e_j)
\bigr)\,
T_f (e_i,\,e_j)
\,dv_g
\nonumber \\
&  &
\ + \
3\int_M
\varphi_k ^2
\|T_f\|^2
\,dv_g\,.
\nonumber
\end{eqnarray}
Using the following formula, 
we replace the curvature term of $N$ 
by the Ricci curvature term of $M$.

\vskip 3ex

\noindent
{\bf Lemma 3}.\ \ 
Let $f$ be a smooth map from $M$ into $N$. 
For any vector fields $X$ and $Y$ on $M$, 
the following equality holds:
\begin{eqnarray}
\lefteqn{
\frac{1}{4}
\nabla_X \nabla_Y \|T_f\|^2
} \label{eq:differential equality} \\
& = &
\sum_{i}
h\bigl( (\nabla_{e_i}\nabla_X df)(Y),\,\sigma_f (e_i)\bigr)
\nonumber \\
&  &
\ + \
\sum_{i,\,j}
h\bigl( (\nabla_X df)(e_i),\,(\nabla_Y df)(e_j)\bigr)\,
T_f (e_i,\,e_j)
\nonumber \\
&  &
\ + \
\sum_{i,\,j}
h\bigl( (\nabla_X df)(e_i),\,df(e_j)\bigr)
h\bigl((\nabla_Y df)(e_i),\,df(e_j)\bigr)
\nonumber \\
&  &
\ + \
\sum_{i,\,j}
h\bigl( (\nabla_X df)(e_i),\,df(e_j)\bigr)
h\bigl(df(e_i),\,(\nabla_Y df)(e_j)\bigr)
\nonumber \\
&  & 
\ -\ 
\frac{2}{m} 
\sum_{i} 
h \bigl( 
(\nabla_{X} df)(e_i)
,\,
df(e_i)
\bigr)
\, 
\sum_{j} 
h \bigl( 
(\nabla_{Y} df)(e_j)
,\,
df(e_j)
\bigr)
\nonumber \\ 
&  &
\ - \
\sum_{i,\,j}
\,
h\bigl({} df(^M\!R(X,\,e_i)(Y)),\,df(e_j)\bigr)\,
T_f (e_i,\,e_j)
\nonumber \\
&  &
\ + \
\sum_{i,\,j}
h\bigl({}^N\!R\bigl(df(X),\,df(e_i)\bigr)df(Y),\,df(e_j)\bigr)\,
T_f (e_i,\,e_j)\,.
\nonumber
\end{eqnarray}
where $\{\,e_i\}$ denotes 
an orthonormal frame on M.

\vskip 3ex

\noindent
{\bf Remark 4}.\
Using integration by parts, we can verify that 
the integral of the first term of the right hand side
in Lemma 3 over $M$ vanishes
for C-stationary maps $f$.

\vskip 3ex

\noindent
{\bf Proof of Lemma 3}.\ \ 
By Lemma 1, We see 
\begin{eqnarray}
\lefteqn{
\frac{1}{4}
\nabla_X \nabla_Y \|T_f\|^2
\ = \
\frac{1}{4}
\nabla_X \nabla_Y 
\left(
\sum_{i,\,j}
T_f(e_i,\,e_j)^2
\right)
}
\label{eq:proof_of_Lemma2_no1} \\
& = &
\sum_{i,\,j}
h\bigl( (\nabla_X\nabla_Y df)(e_i),\,df(e_j)\bigr)
T_f(e_i,\,e_j)
\nonumber \\
&  &
\ + \
\sum_{i,\,j}
h\bigl( (\nabla_Y df)(e_i),\,(\nabla_X df)(e_j)\bigr)
T_f(e_i,\,e_j)
\nonumber \\
&  &
\ + \
\sum_{i,\,j}
h\bigl( (\nabla_Y df)(e_i),\,df(e_j)\bigr)
h\bigl((\nabla_X df)(e_i),\,df(e_j)\bigr)
\nonumber \\
&  &
\ + \
\sum_{i,\,j}
h\bigl( (\nabla_Y df)(e_i),\,df(e_j)\bigr)
h\bigl(df(e_i),\,(\nabla_X df)(e_j)\bigr) 
\nonumber \\ 
&  &
\ - \ 
\frac{2}{m}
\sum_{i}
h\bigl( (\nabla_Y df)(e_i),\,df(e_i)\bigr)
\sum_{j}
h\bigl((\nabla_X df)(e_j),\,df(e_j)\bigr)\,.
\nonumber 
\end{eqnarray}
By the Ricci formula, we have 
\begin{eqnarray}
\lefteqn{
(\nabla_X \nabla_Y df)(e_i)
\ = \
(\nabla_X \nabla_{e_i} df)(Y)
} \label{eq:proof_of_Lemma2_no2}  \\
& = &
(\nabla_{e_i}\nabla_X df)(Y)
\ - \
\,df( {}^M\!R(X,\,e_i)(Y))
\ + \
{}^N\!R\bigl(df(X),\,df(e_i)\bigr)df(Y)
\nonumber 
\end{eqnarray}
Furthermore by Lemma 2 we have 
\begin{eqnarray}
\qquad 
\sum_{i,\,j}
h\bigl( (\nabla_{e_i}\nabla_X df)(Y),\,df(e_j)\bigr)\,
T_f (e_i,\,e_j)
& = &
\sum_{i}
h\bigl( (\nabla_{e_i}\nabla_X df)(Y),\,\sigma_f (e_i)\bigr). 
\label{eq:proof_of_Lemma2_no3} 
\end{eqnarray}
By $(\ref{eq:proof_of_Lemma2_no1})$, $(\ref{eq:proof_of_Lemma2_no2})$ and $(\ref{eq:proof_of_Lemma2_no3})$, 
we have the equality in Lemma 3.
$\square$

\vskip 3ex

%
By Lemma 3 with $X$ $=$ $Y$ $=$ $Z_k$, we have 
\begin{eqnarray*}
& & 
\sum_{i,\,j}
h\bigl( (\nabla_{Z_k} df)(e_i),\,(\nabla_{Z_k} df)(e_j)\bigr)\,
T_f (e_i,\,e_j)
\nonumber \\
&  &
\ + \
\sum_{i,\,j}
h\bigl( (\nabla_{Z_k} df)(e_i),\,df(e_j)\bigr)
h\bigl((\nabla_{Z_k} df)(e_i),\,df(e_j)\bigr)
\nonumber \\
&  &
\ + \
\sum_{i,\,j}
h\bigl( (\nabla_{Z_k} df)(e_i),\,df(e_j)\bigr)
h\bigl(df(e_i),\,(\nabla_{Z_k} df)(e_j)\bigr)
\nonumber \\
&  & 
\ - \ 
\frac{2}{m} 
\sum_{i} 
h \bigl( 
(\nabla_{Z_k} df)(e_i)
,\,
df(e_i)
\bigr)
\, 
\sum_{j} 
h \bigl( 
(\nabla_{Z_k} df)(e_j)
,\,
df(e_j)
\bigr)
\nonumber \\ 
&  &
\ + \
\sum_{i,\,j}
h\bigl({}^N\!R\bigl(df(Z_k),\,df(e_i)\bigr)df(Z_k),\,df(e_j)\bigr)\,
T_f (e_i,\,e_j)\,.
\nonumber \\
& = &
\frac{1}{4}
\nabla_{Z_k} \nabla_{Z_k} \|T_f\|^2
\nonumber \\ 
&  & 
-\,\sum_{i}
h\bigl( (\nabla_{e_i}\nabla_{Z_k} df)(Z_k),\,\sigma_f (e_i)\bigr)
\nonumber \\
&  &
\ + \
\sum_{i,\,j}
\,
h\bigl({} df(^{\mathrm{S}^m}\!R(Z_k,\,e_i)(Z_k)),\,df(e_j)\bigr)\,
T_f (e_i,\,e_j)
\nonumber 
\end{eqnarray*}
Then by the equalities $(\ref{eq:stability___02})$ 
\begin{eqnarray}
0 & \leq &
\sum_{k=1}^{m+1}
{\mathrm L}(df(Z_k),\,df(Z_k))
\label{eq:stability___03}
\\
& = &
\frac{1}{4}
\int_M
\sum_{k=1}^{m+1}
\nabla_{Z_k}\nabla_{Z_k} \|T_f\|^2 \,dv_g
\nonumber \\
& &
\ - \
\int_M
\sum_{k=1}^{m+1}
\sum_{i}
h \bigl(
(\nabla_{e_i} \nabla_{Z_k}df)(Z_k)
,\,
\sigma_f(e_i)
\bigr)
\,dv_g
\nonumber \\
&  &
\ + \
\int_M
\sum_{k=1}^{m+1}
\sum_{i,\,j}
h( df({}^{S^m}\!R (Z_k,\,e_i)
Z_k),\,df(e_j)\bigr)\,
T_f (e_i,\,e_j)
dv_g
\nonumber \\
&  &
\ - \
6
\int_M
\sum_{k=1}^{m+1}
\varphi_k
\sum_{i,\,j}
h \bigl(
(\nabla_{e_i}df)(Z_k)
,\,
df(e_j)
\bigr)\,
T_f (e_i,\,e_j)
\,dv_g
\nonumber \\
&  &
\ + \
3\,
\int_M
\sum_{k=1}^{m+1}
\varphi_k ^2\,
\|T_f\|^2
\,dv_g.
\nonumber 
\end{eqnarray}
We write terms on the right hand side as I, II, III, IV and V 
respectively.
To calculate these terms, 
we use the following lemma. 

\vskip 3ex

\noindent 
{\bf Lemma 4}.\\ \\ 
\qquad (a)\ \ 
${\displaystyle 
\sum_{k=1}^{m+1}
\nabla_{Z_k}
\nabla_{Z_k}
\ = \
\bigtriangleup
}$,\ 
where $\bigtriangleup$ denotes the Laplacian on $M$ $=$ $\mathrm{S}^m$. 
\\ 
\qquad (b)\ \ 
${\displaystyle
\sum_{k=1}^{m+1}
g(e_i,\,Z_k) Z_k
}$ $=$ $e_i$. 
\\ 
\qquad (c)\ \ 
${\displaystyle
\sum_{k=1}^{m+1}
\sum_{i}
h\bigl(df(Z_k),\,df(e_i)\bigr)\,
T_f (Z_k,\,e_i)
}$ $=$ $\|T_f\|^2$.

\vskip 3ex

\noindent 
{\bf Proof of Lemma 4}.\ \ 
We first prove that 
${\displaystyle 
\sum_{k=1}^{m+1}
\nabla_{Z_k}
\nabla_{Z_k}
}$  
does not depend on the choice of $E_k$. 
Take any two sets of orthonormal parallel vector fields 
$\{E_k\}$, $\{\overline{E}_k \}$ on $\mathbb{R}^{m+1}$\,. 
Then there exists an orthogonal matrix $(a_{kp})$
such that
${\displaystyle 
E_k =  \sum_{p = 1}^{m+1}
a_{kp}\ \overline{E}_p
}$.
%
Hence
${\displaystyle 
Z_k =  \sum_{p = 1}^{m+1}
a_{kp}\ \overline{Z}_k
}$, where $\overline{Z}_k = p(\overline{E}_k)$. 
Thus 
we have
\begin{eqnarray*}
\sum_{k=1}^{m+1}
\nabla_{Z_k}
\nabla_{Z_k}
& = &
\sum_{p=1}^{m+1}\sum_{q=1}^{m+1}
\sum_{k=1}^{m+1} a_{kp} a_{kq}
\nabla_{\overline{Z}_p}
\nabla_{\overline{Z}_q}
\\
& = &
\sum_{\ell =1}^{m+1}
\nabla_{\overline{Z}_\ell}
\nabla_{\overline{Z}_\ell}\,,
\end{eqnarray*}
which implies that 
this operator does not depend on the choice of $E_k$. 
At a point $x$ on $\mathrm{S}^m$,
take orthonormal parallel vector fields 
$E_1$, $\cdots$ ,\,$E_{m+1}$ on $\mathbb{R}^{m+1}$
such that $Z_k(x)$ $(k=1,\,\cdots\,,\,m)$ is
an orthonormal basis of $T_x \mathrm{S}^m$ and
$Z_{m+1}(x) = 0$.
Then the equality in (a) holds clearly. 

To prove (b) and (c),
we have only to prove 
${\displaystyle 
\sum_{k=1}^{m+1}
g(X,\,Z_k) Z_k
}$ and 
${\displaystyle
\sum_{k=1}^{m+1}
\sum_{i}
h\bigl(df(Z_k),\,df(e_i)\bigr)\,
T_f (Z_k,\,e_i)
}$ do not depend on the choice of $E_k$, 
which we can verify easily. 
$\square$

\vskip 5ex

By Lemma 4 (a), we have 
\begin{eqnarray*}
\mathrm{I}
& = &
\frac{1}{4}
\int_M
\sum_{k=1}^{m+1}
\nabla_{Z_k}
\nabla_{Z_k}
\|T_f\|\,dv_g
\ {=} \ 
\frac{1}{4}
\int_{M}
\bigtriangleup
\|T_f\|\,dv_g
\ = \
0\,.
\end{eqnarray*}
Using the integration by parts, we have 
\begin{eqnarray*}
\mathrm{II} & = & 
\int_M 
\sum_{i}
h \bigl(
(\nabla_{e_i} \nabla_{Z_k}df)(Z_k)
,\,
\sigma_f(e_i)
\bigr)\,dv_g
\nonumber \\ 
& = & 
-\,
\int_M 
h \bigl(
(\nabla_{Z_k}df)(Z_k)
,\,
\mathrm{div}_g \sigma_f 
\bigr)\,dv_g
\nonumber \\ 
& = & 0\,, 
\end{eqnarray*}
since $f$ is a C-stationary map, i.e., 
$\mathrm{div}_g\sigma_f$ $=$ $0$.

Considering the fact that
${}^{\mathrm{S}^m}\!R(U,\, V)W = g(V,W)U - g(U,W)V$,
we get by Lemma 4 (b)
\begin{eqnarray*}
\mathrm{III}
& = &
\int_M 
\sum_{k=1}^{m+1}
\sum_{i,\,j}
h \Bigl(
df\bigl({}^{\mathrm{S}^m}\!R(Z_k, e_i)Z_k\bigr),\,df(e_j)
\Bigr)\,T_f(e_i,\,e_j)
\\
& = &
\int_M
\sum_{i,\,j}
\sum_{k=1}^{m+1}
h \Bigl(
df\bigl(g(e_i, Z_k)Z_k - g(Z_k, Z_k)e_i\bigr),\,df(e_j)
\Bigr)
\,T_f(e_i,\,e_j)
\\
& = &
-\,(m-1)
\int_M
\sum_{i,\,j}
h \bigl(
df(e_i)
,\,
df(e_j)
\bigr)\,
T_f (e_i,\,e_j)
\\
& = &
-\,(m-1)
\int_M
\|T_f\|^2\,.
\end{eqnarray*}
We used here 
\begin{eqnarray*}
\sum_{k=1}^{m+1} g(e_i,\,Z_k)Z_k
& = & 
p\left(
\sum_{k=1}^{m+1} g(e_i,\,E_k)E_k
\right)
\ = \ p(e_i)\ = \ e_i\,.
\end{eqnarray*} 
For the calculation of the term IV, let us define
\begin{eqnarray*}
\gamma_k(X)
\:\,=\
h\bigl(df(Z_k),\,\sigma_f(X)\bigr).
\end{eqnarray*}
Then by $(\ref{eq:equality_for_vector_field_k})$ and Lemma 2 we have
\begin{eqnarray}
\lefteqn{
\sum_{i,\,j}
h \bigl(
(\nabla_{e_i}df)(Z_k)
,\,
df(e_j)
\bigr)\,
T_f (e_i,\,e_j)
} \nonumber \\
& = &
\sum_{i}
h\bigl( (\nabla_{e_i}df)(Z_k)
,\,
\sigma_f(e_i)
\bigr)
\nonumber \\
& = &
\sum_{i}
\Bigl\{
h\bigl(
 \nabla_{e_i}\bigl(df(Z_k)\bigr)
,\,
\sigma_f(e_i)
\bigr)
\ + \
\varphi_k
\sum_{i}
h\bigl(
df(e_i)
,\,
\sigma_f(e_i)
\bigr)
\Bigr\}
\nonumber \\
&=&  \sum_{i}
(\nabla_{e_i}\gamma_k)(e_i)
\ - \
h\bigl(
df(Z_k)
,\,
\sum_{i}
(\nabla_{e_i}\sigma_f)(e_i)
\bigr)
\ + \
\varphi_k
\sum_{i}
h\bigl(
df(e_i)
,\,
\sigma_f(e_i)
\bigr)
\nonumber \\
& = &
\mathrm{div}\,\gamma_k
\ + \
\varphi_k
\|T_f\|^2, 
\nonumber
\end{eqnarray}
since ${\displaystyle
\sum_{i}\bigl(\nabla_{e_i}\sigma_f\bigr)(e_i)
\ = \
\mathrm{div} \sigma_f
\ = \ 0
}$. 
Hence we have
\begin{eqnarray}
\lefteqn{
\sum_{k=1}^{m+1}
\varphi_k
\sum_{i,\,j}
h \bigl(
(\nabla_{e_i}df)(Z_k)
,\,
df(e_j)
\bigr)\,
T_f (e_i,\,e_j)
} \label{eq:term_00} \\
& = &
\sum_{k=1}^{m+1}
\varphi_k
\mathrm{div}\,\gamma_k
\ + \
\sum_{k=1}^{m+1}
\varphi_k^2
\|T_f\|^2 
\nonumber \\
& = &
\mathrm{div}\,\left(
\sum_{k=1}^{m+1}
\varphi_k \gamma_k
\right)
\ - \
\sum_{k=1}^{m+1} \sum_i
e_i(\varphi_k)\,\gamma_k(e_i)\,
\ + \
\sum_{k=1}^{m+1}
\varphi_k^2
\|T_f\|^2\,.
\nonumber
\end{eqnarray}
We see 
\begin{eqnarray*}
\sum_i e_i(\varphi_k)\,e_i
& = & 
\sum_i e_i \left(\,\langle E_k,\,\nu \rangle\,\right) \,e_i
\ = \ 
\sum_i\langle E_k,\,{}^{\mathbb{R}^{n+1}}\nabla_{e_i}\nu \rangle \,e_i
\ = \ 
\sum_i \langle E_k,\,e_i \rangle \,e_i
\ = \ 
p(E_k)
\ = \ 
Z_k\,, 
\nonumber 
\end{eqnarray*}
since $E$ is parallel and 
${}^{\mathbb{R}^{n+1}}\nabla_{e_i}\nu = e_i$ 
where ${}^{\mathbb{R}^{n+1}}\nabla$ denotes the standard connection on $\mathbb{R}^{n+1}$. 
Hence by Lemma 4 (c), we have 
\begin{eqnarray}
\lefteqn{
\sum_{k=1}^{m+1} \sum_i
e_i(\varphi_k)\,\gamma_k(e_i)\,
\ = \
\sum_{k=1}^{m+1}
\gamma_k
\bigl(
\sum_i
e_i(\varphi_k)\,e_i
\bigr)
} \label{eq:varphi_01} \\
& = &
\sum_{k=1}^{m+1}
\gamma_k
\bigl(
Z_k
\bigr)
\ = \
\sum_{k=1}^{m+1}
h \bigl(
df(
Z_k
)
,\,
\sigma_f(
Z_k
)
\bigr)
\nonumber \\
& = &
\sum_{k=1}^{m+1}
\sum_i
h \bigl(
df(
Z_k
)
,\,
df(e_i)
\bigr)\,
T_f (Z_k,\,e_j)
\ = \
\|T_f\|^2\,.
\nonumber
\end{eqnarray}
We see 
\begin{eqnarray}
\sum_{k=1}^{m+1}
\varphi_k^2
= 
\sum_{k=1}^{m+1}
g(E_k,\,\nu )^2
\label{eq:varphi_02} 
\ = \ 
g\bigl(\nu,\,\sum_{k=1}^{m+1} g(\nu,\,E_k) E_k \bigr)
\ = \ 
g(\nu,\,\nu )
\ = \ 1\,.
\end{eqnarray}
Then by
$(\ref{eq:term_00})$, $(\ref{eq:varphi_01})$ 
and $(\ref{eq:varphi_02})$,
we obtain $\mathrm{IV} = 0$\,.
By $(\ref{eq:varphi_02})$, we have 
\begin{eqnarray*}
V
& = &
3 \int_M \|T_f\|^2\,dv_g\,.
\end{eqnarray*}

Finally substituting I through V into $(\ref{eq:stability___03})$, 
we have 
\begin{eqnarray*}
0 & \leq &
\sum_{k=1}^{m+1}
\mathrm{L}\bigl (df(Z_k),\,df(Z_k)\bigr)
\\
& = &
-\,(m-1)\int_M \|T_f\|^2\,dv_g
\ + \
3 \int_M \|T_f\|^2\,dv_g
\\
& = &
(4-m) \int_M \|T_f\|^2\,dv_g,
\end{eqnarray*}
i.e.,
\begin{eqnarray*}
(m - 4) \int_M \|T_f||^2\,dv_g\ \leq 0\,.
\end{eqnarray*}
Hence, if $m \geq 5$,
then we get $\|T_f\| = 0$,
i.e. $f$ is a weakly conformal map.
$\square$

\section{Stable maps into spheres}

In this section we prove Theorem 2.

\vskip 5ex

\noindent
{\bf Proof of Theorem 2}.\ \
We use notations similar to those in the proof of Theorem 1.
Since the standard sphere $\mathbb{S}^n$ is
a submanifold of the Euclidean space $\mathbb{R}^{n+1}$,
we may consider that
the tangent spaces of standard sphere $\mathbb{S}^n$
at $y$ $\in$ $\mathbb{S}^n$
is a subspace of the linear space $\mathbb{R}^{n+1}$ $\simeq$
$T_y\mathbb{R}^{n+1}$.
Let $p=p_y$ denotes the canonical projection
from $\mathbb{R}^{n+1}$
onto the tangent space $\mathrm{T}_y\mathbb{S}^{n}$.

Let $E$ be normal parallel vector field on $\mathbb{R}^{n+1}$, 
and let $Z$ be the vector field which is
the image by the projection $p$ of
$E$.
Then we define a smooth section $W$ of the pull-back bundle 
$f^{-1}T\mathbb{S}^n$ by
\begin{eqnarray*}
W_x
\ :\,=\ Z_{f(x)}
\ = \
p_{f(x)}\bigl(E_{f(x)}\bigr)
\end{eqnarray*} 
for $x$ $\in$ $M$.
Then we see
\begin{eqnarray}
\nabla_{e_i}W
& = &
{}^{f^{-1}T\mathbb{S}^n}\nabla_{e_i}W
\ = \
{}^{\mathbb{S}^n}\nabla_{df(e_i)}Z
\ = \
-\,\varphi\,df(e_i)\,,
\label{eq:varphi_03} 
\end{eqnarray}
where 
\begin{eqnarray}
\varphi 
\ = \
h(E,\,\nu )
\nonumber
\end{eqnarray}
and $\nu = \nu_x$ is the unit outer normal vector at $f(x)$. 
(Note that $\nu_x = f(x)$ in the case of $N$ $=$ $\mathrm{S}^n$. )

For simplicity, 
we use the notation $Z$ instead of $W$, 
and then by $(\ref{eq:varphi_03})$ we have 
\begin{eqnarray}
\nabla_{e_i}Z
\ = \
-h(E,\,\nu )\,df(e_i)\,.
\label{eq:varphi_04} 
\end{eqnarray}

\vskip 3ex

%
Then by $(\ref{eq:varphi_04})$, 
we have
\begin{eqnarray}
\mathrm{L}(Z,\,Z)
& = &
\int_M
\sum_{i,\,j}
h \bigl(
\nabla_{e_i}Z
,\,
\nabla_{e_j}Z
\bigr)
\,
T_f (e_i,\,e_j)
\,dv_g
\label{eq:2nd_variation_to_sphere} \\
&  &
\ + \
\int_M
\sum_{i,\,j}
h \bigl(
\nabla_{e_i}Z
,\,
df(e_j)
\bigr)^2
\,dv_g
\nonumber \\
&  &
\ + \
\int_M
\sum_{i,\,j}
h \bigl(
\nabla_{e_i}Z
,\,
df(e_j)
\bigr)
\,
h \bigl(
df(e_i)
,\,
\nabla_{e_j}Z
\bigr)
\,dv_g
\nonumber \\
&  & 
\ - \ 
\frac{2}{m}
\int_M 
\left\{
\sum _{i} 
h \bigl( 
\nabla_{e_i}Z
,\,
df(e_i)
\bigr)
\right\}^2
\,dv_g
\nonumber \\
&  &
\ + \
\int_M
\sum_{i,\,j}
h\bigl({}^{\mathrm{S}^n}\!
R\left(
Z, \,
df(e_i)
\right)
Z,\,
df(e_j)
\bigr)\,
T_f (e_i,\,e_j)\,
dv_g
\nonumber \\
& = &
3\int_M
h(E,\,\nu )^2
\sum_{i,\,j}
h \bigl(df(e_i),\,df(e_j)\bigr)\,
T_f (e_i,\,e_j)
\,dv_g
\nonumber \\
&  &
\ - \
\int_M
h(Z,\,Z )
\sum_{i,\,j}
h \bigl(
df(e_i)
,\,
df(e_j)
\bigr)\,
T_f (e_i,\,e_j)
\,dv_g
\nonumber \\
&  &
\ + \
\int_M
\sum_{i,\,j}
h \bigl(
Z
,\,
df(e_i)
\bigr)
h \bigl(
Z
,\,
df(e_j)
\bigr)\,
T_f (e_i,\,e_j)
\,dv_g
\nonumber \\
& = &
3 \int_M
h(E,\,\nu )^2
\|T_f\|^2 \,dv_g
\nonumber \\
&  &
\ - \
\int_M
h(Z,\,Z )
\|T_f\|^2
\,dv_g
\nonumber \\
&  &
\ + \
\int_M
\sum_{i,\,j}
h \bigl(
Z
,\,
df(e_i)
\bigr)
h \bigl(
Z
,\,
df(e_j)
\bigr)\,
T_f (e_i,\,e_j)
\,dv_g\,.
\nonumber
\end{eqnarray}
Take orthonormal parallel vector fields 
$E_1$,\,$\cdots$\,,\,$E_{n+1}$ 
on $\mathbb{R}^{n+1}$\,. 
The stability of the C-stationary map $f$ 
implies the inequality 
\begin{eqnarray*}
L(Z_k,\,Z_k ) \geq 0,
\ \text{hence}\ 
\sum_{k=1}^{n+1} L(Z_k,\,Z_k ) \geq 0\,. 
\end{eqnarray*}
Then by $(\ref{eq:2nd_variation_to_sphere})$, 
we have 
\begin{eqnarray}
\lefteqn{
0 \leq 
\sum_{k=1}^{n+1}
\mathrm{L}(Z_k,\,Z_k)
} 
\label{eq:stability___05}  \\
& = &
3 \int_M
\sum_{k=1}^{n+1}
h(E_k,\,\nu )^2
\|T_f\|^2 \,dv_g
\nonumber \\
&  &
\ - \
\int_M
\sum_{k=1}^{n+1}
h(Z_k,\,Z_k )
\|T_f\|^2
\,dv_g
\nonumber \\
&  &
\ + \
\int_M
\sum_{i,\,j}
\sum_{k=1}^{n+1}
h \bigl(
Z_k
,\,
df(e_i)
\bigr)\,
h \bigl(
Z_k
,\,
df(e_j)
\bigr)\,
T_f (e_i,\,e_j)
%
%
%
\,dv_g\,.
\nonumber
\end{eqnarray}

\vskip 3ex 

\noindent 
To calculate $L(Z_k,\,Z_k )$, 
we first give the following lemma.

\vskip 3ex

\noindent 
{\bf Lemma 5}. \qquad 
${\displaystyle 
\sum_{k=1}^{n+1}
h(
Z_k,\,
Z_k
)
\ = \
n
}$

\vskip 3ex

\noindent 
{\bf Proof of Lemma 5}.\ \ 
We first prove that 
${\displaystyle 
\sum_{k=1}^{n+1}
h(
Z_k,\,
Z_k
)
}$ does not depend on the choice of $E_k$. 
Take any two sets of orthonormal parallel vector fields $\{E_k\}$, 
$\{\overline{E}_k \}$ on $\mathbb{R}^{n+1}$\,. 
Then there exists an orthogonal matrix $(a_{kp})$
such that 
${\displaystyle 
Z_k =  \sum_{p = 1}^{n+1}
a_{kp}\ \overline{Z}_p
}$. 
Thus we have
\begin{eqnarray*}
\sum_{k=1}^{n+1}
h(
Z_k,\,
Z_k
)
& = &
\sum_{p=1}^{n+1}\sum_{q=1}^{n+1}
\sum_{k=1}^{n+1} a_{kp} a_{kq}
h(
\overline{Z}_p,\,
\overline{Z}_q
)
\ = \
\sum_{\ell =1}^{n+1}
h(
\overline{Z}_{\ell},\,
\overline{Z}_{\ell}
)
\end{eqnarray*}
which implies that 
this quantity does not depend on the choice of $E_k$. 
At a point $x$ on $\mathrm{S}^m$, 
take orthonormal parallel frame $\{E_1,\,\cdots\,,\,E_{n+1}\}$ on $\mathbb{R}^{n+1}$ 
such that $\{Z_1(x),\,\cdots\,,\,Z_m(x)\}$ is
an orthonormal base of $T_x \mathrm{S}^m$ and
$Z_{n+1}(x) = 0$, and then we have Lemma 5. 
$\square$

\vskip 5ex

Since 
$Z_k$ is the projection of $E_k$ onto
the tangent space of the sphere $\mathbb{S}^n$,
we have
\begin{eqnarray*}
\sum_{k=1}^{n+1} h(E_k,\,\nu )^2
& = &
h(\nu,\,\nu )
\ = \
1 
\end{eqnarray*}
and 
\begin{eqnarray*}
\lefteqn{
\sum_{k=1}^{n+1}
\sum_{i,\,j}
h \bigl(
Z_k
,\,
df(e_i)
\bigr)\,
h \bigl(
Z_k
,\,
df(e_j)
\bigr)\,
T_f (e_i,\,e_j)
}  \\
& = &
\sum_{i,\,j}
h\bigl (df(e_i),\,
\sum_{k=1}^{n+1}
h(df(e_j),\,Z_k)\,
Z_k\bigr)\,
T_f (e_i,\,e_j)
\\
& = &
\sum_{i,\,j}
h \bigl(df(e_i),\,df(e_j)\bigr)\,
T_f (e_i,\,e_j)
\ = \
\|T_f\|^2\,.
\end{eqnarray*}
Therefore we get 
by Lemma 5 and $(\ref{eq:stability___05})$ 
\begin{eqnarray*}
0 & \leq &
\sum_{k=1}^{n+1} \mathrm{L}(Z_k,\,Z_k)
\\
& = &
3 \int_M \|T_f\|^2\,dv_g
\ - \
n \int_M \|T_f\|^2\,dv_g
\ + \
\int_M \|T_f\|^2\,dv_g
\\
& = &
(4-n) \int_M \|T_f\|^2\,dv_g\,.
\end{eqnarray*}
Thus we obtain
\begin{eqnarray*}
(n - 4) \int_M \|T_f\|^2\ \leq 0\,.
\end{eqnarray*}
Hence, if $n \geq 5$,
then we get $\|T_f\| = 0$\,., 
i.e. $f$ is a weakly conformal map.
$\square$

\vskip 5ex

\noindent
Faculty of Culture and Education \\
Saga University \\
Saga 840 \\
Japan

\vskip 5ex

\noindent
Faculty of Science \\
Yamaguchi University \\
Yamaguchi 753-8512 \\
Japan

\newpage


\noindent
{\bf Appendix} 

\vskip 5ex

In this appendix, 
we give, 
for reader's convenience, 
proofs of the first variation formula (Proposition 1) 
and the second variation formula (Proposition 2). 

\vskip 3ex  

\noindent 
{\bf Proof of Proposition 1}. 
We calculate 
$\frac{\partial}{\partial t}\left\|f^*_t h\right\|^2$ at any fixed point
$x_0 \in M$.
The connection $\nabla$ is trivially extended to
a connection on $(-\varepsilon,\,\varepsilon ) \times M$.
The frame $e_i$ is also trivially extended to a frame
on $(-\varepsilon,\,\varepsilon ) \times $ (the domain of the frame),
and we use the same notation $e_i$.
Then we see $\nabla_{e_i}\frac{\partial}{\partial t}$ $=$
$\nabla_{\frac{\partial}{\partial t}}e_i$ $=$ $0$
on $(-\varepsilon,\,\varepsilon ) \times M$\,.
Using a normal coordinate at $x_0$, 
We can assume $\nabla_{e_i}{e_j} = 0$  at $x_0$ for any $i$, $j$\,.
Since $(dF)_{(t,\,x)}\bigl((e_i)_{(t,\,x)}\bigr)$
$=$ $(dF)_x\bigl((e_i)_x\bigr)$, 
we denote them by $dF(e_i)$ simply. 
Note 
\begin{eqnarray}
\nabla_{{\textstyle \frac{\partial}{\partial t}}} \bigl(dF (e_i)\bigr)
\ = \
\nabla_{e_i}
\left(dF \bigl({\textstyle \frac{\partial}{\partial t}}\bigr) \right)\,,
\label{eqn:commutative}
\end{eqnarray}
since $\left[{\textstyle \frac{\partial}{\partial t}},\,e_i\right]$ $=$ $0$.
Then using Lemma 1 (a) and (d), we have 
\begin{eqnarray*}
\lefteqn{
\frac{\partial}{\partial t}
\|T_{f_t}\|^2 
\ = \ 
\frac{\partial}{\partial t}
\sum_{i,\,j} T_F(e_i,\,e_j)^2 
} \\ 
& = & 
2\,
\sum_{i,\,j} 
\frac{\partial T_F(e_i,\,e_j) }{\partial t}
T_F(e_i,\,e_j)
\\ 
& = & 
2\,
\sum_{i,\,j} 
\left\{
\frac{\partial}{\partial t}
h\bigl(dF(e_i),\,dF(e_j)\bigr)
-
\frac{1}{m} 
\frac{\partial \|dF\|^2\,}{\partial t}
g(e_i,\,e_j) 
\right\}
T_F(e_i,\,e_j)
\\ 
& = & 
2\,
\sum_{i,\,j} 
\left\{
\frac{\partial}{\partial t}
h\bigl(dF(e_i),\,dF(e_j)\bigr)
\right\}
T_F(e_i,\,e_j)
\ - \ 
\frac{2}{m} 
\frac{\partial \|dF\|^2\,}{\partial t}
\sum_{i,\,j} 
g(e_i,\,e_j)\,
T_F(e_i,\,e_j)
\\ 
& = & 
2\,
\sum_{i,\,j} 
\left\{
\frac{\partial}{\partial t}
h\bigl(dF(e_i),\,dF(e_j)\bigr)
\right\}
T_F(e_i,\,e_j)
\\ 
& = & 
4\,
\sum_{i,\,j} 
h\bigl({\textstyle \nabla_{{\textstyle \frac{\partial}{\partial t}}}}\bigl(dF(e_i)\bigr),\,dF(e_j)\bigr)
T_F(e_i,\,e_j)
\\ 
& = & 
4\,
\sum_{i,\,j} 
h\bigl(\nabla_{e_i}\bigl(dF({\textstyle \frac{\partial}{\partial t}})\bigr),\,
dF(e_j)\bigr)
T_F(e_i,\,e_j)
\\ 
& = & 
4\,
\sum_{i} 
h\bigl(\nabla_{e_i}\bigl(dF({\textstyle \frac{\partial}{\partial t}})\bigr),\,
\sigma_{F} (e_i) \bigr).
\end{eqnarray*}
The last equality follows from 
Lemma 2 
for $Z$ $=$ $\nabla_{e_i} 
\left(dF \bigl({\textstyle \frac{\partial}{\partial t}}\bigr) \right)$ 
and for $W$ $=$ $e_i$\,. 
Integrate it over $M$ and let $t=0$.
Then
using integration by parts,
we obtain the first variation formula. $\square$

\vskip 3ex

\noindent 
{\bf Proof of Proposition 2}.\ \ 
We calculate 
$\frac{\partial^2}{\partial s \partial t}\left\|f^*_{s,\,t} h\right\|^2$ at any fixed point
$x_0 \in M$.
The connection $\nabla$ is trivially extended to
a connection on $(-\varepsilon,\,\varepsilon ) \times (-\delta,\,\delta )
\times M$.
We use the same notation $\nabla$ for this connections.
The frame $e_i$ is also trivially extended to a frame
on $(-\varepsilon,\,\varepsilon ) \times (-\delta,\,\delta ) \times $
(the domain of the frame), denoted by the same notation $e_i$.
Then we see
\begin{eqnarray*}
\nabla_{\textstyle \frac{\partial}{\partial s}}e_i
& = &
\nabla_{e_i}\frac{\textstyle \partial}{\partial s}
\ = \ 0,
\\
\nabla_{\textstyle \frac{\partial}{\partial t}}e_i
& = &
\nabla_{e_i}{\textstyle \frac{\partial}{\partial t}}
\ = \ 0,
\\
\nabla_{{\textstyle \frac{\partial}{\partial s}}}{\textstyle \frac{\partial}
{\partial t}}
& = &
\nabla_{{\textstyle \frac{\partial}{\partial t}}}{\textstyle \frac{\partial}
{\partial s}}
\ = \ 0
\end{eqnarray*}
%
on $(-\varepsilon,\,\varepsilon ) \times (-\delta,\,\delta ) \times M$\,.
Take and fix any point $x_0$ $\in$ $M$.
We
calculate
$\frac{\partial}{\partial t} \|f_t^*h\|^2$
at $x$ $=$ $x_0$. 
Since
we can assume $\nabla_{e_i}e_j = 0$
at $x_0$
for any $i$, $j$\,,
we see at $x_0$,
\begin{eqnarray*}
\nabla_{{\textstyle \frac{\partial}{\partial s}}} \bigl(dF (e_i)\bigr)
& = &
\nabla_{e_i} \left(dF \bigl({\textstyle \frac{\partial}{\partial s}}\bigr)
\right)\,,
\label{eq:vect_1}
\\
\nabla_{{\textstyle \frac{\partial}{\partial t}}} \bigl(dF (e_i)\bigr)
& = &
\nabla_{e_i} \left(dF \bigl({\textstyle
\frac{\partial}{\partial t}}\bigr) \right)\,.
\label{eq:vect_2}
\end{eqnarray*}
Then we have 
\begin{eqnarray}
&  & 
\frac{1}{4} 
\frac{\partial^2}{\partial s \partial t}
\|T_{f_{s,\,t}}\|^2 
\label{eq:second_variation} \\ 
& = & 
\frac{1}{4} 
\frac{\partial^2}{\partial s \partial t}
\sum_{i,\,j} 
T_F(e_i,\,e_j)^2
\nonumber \\ 
& = & 
\frac{1}{2} 
\,\sum_{i,\,j} 
\left\{
\frac{\partial^2 T_F(e_i,\,e_j)}{\partial s\partial t}
T_{F}(e_i,\,e_j)
\right\}
\ + \ 
\frac{1}{2} 
\,\sum_{i,\,j} 
\frac{\partial T_F(e_i,\,e_j)}{\partial s}
\frac{\partial T_F(e_i,\,e_j)}{\partial t}
\nonumber \\ 
& \stackrel{\scriptscriptstyle {\text def}}{=} & 
\mathrm{I}_1\ + \ \mathrm{I}_2
\nonumber  
\end{eqnarray}
We have 
\begin{eqnarray}
& & 
\label{eq:I_1_1} \\ 
\mathrm{I}_1
& = & 
\frac{1}{2}\sum_{i,\,j}
\frac{\partial^2}{\partial s \partial t}
\left\{
h\bigl(dF(e_i),\,dF(e_j)\bigr)
\ - \ 
\frac{1}{m} 
\|dF\|^2\,g(e_i,\,e_j)
\right\}\,
T_F(e_i,\,e_j)
\nonumber \\ 
& = & 
\frac{1}{2}\sum_{i,\,j}
\left\{
\frac{\partial^2}{\partial s \partial t}
h\bigl(dF(e_i),\,dF(e_j)\bigr)
\right\}\,
T_F(e_i,\,e_j)
\ - \ 
\frac{1}{2m} 
\frac{\partial^2\|dF\|^2}{\partial s \partial t}
\sum_{i,\,j}
g(e_i,\,e_j) 
T_F(e_i,\,e_j)
\nonumber \\ 
& = & 
\frac{1}{2}\sum_{i,\,j}
\left\{
\frac{\partial^2}{\partial s \partial t}
h\bigl(dF(e_i),\,dF(e_j)\bigr)
\right\}\,
T_F(e_i,\,e_j)
\qquad
\left(\ \text{by\ Lemma\ 1\,(d)}\ \right)
\nonumber \\ 
& = & 
\frac{1}{2}\sum_{i,\,j}
\left\{
\frac{\partial^2}{\partial s \partial t}
h\bigl(dF(e_i),\,dF(e_j)\bigr)
\right\}\,
T_F(e_i,\,e_j)
\nonumber \\ 
& = & 
\sum_{i,\,j}
\left\{
h\bigl(\nabla_{{\textstyle \frac{\partial}{\partial s}}}
\nabla_{{\textstyle \frac{\partial}{\partial t}}}\bigl(dF(e_i)\bigr),\,
dF(e_j)\bigr)
\right\}\,
T_F(e_i,\,e_j)
\nonumber \\ 
&  & 
\qquad\qquad
\ + \ 
\sum_{i,\,j}
\left\{
h\bigl(\nabla_{{\textstyle \frac{\partial}{\partial s}}}\bigl(dF(e_i)\bigr),\,
\nabla_{{\textstyle \frac{\partial}{\partial t}}}\bigl(dF(e_j)\bigr)\,\bigr)
\right\}\,
T_F(e_i,\,e_j)
\nonumber 
\end{eqnarray}
We get 
\begin{eqnarray}
\nabla_{\textstyle \frac{\partial}{\partial s}}
\nabla_{\textstyle \frac{\partial}{\partial t}}
\bigl(dF(e_i)\bigr)
& = &
\bigl(
\nabla_{\textstyle \frac{\partial}{\partial s}}
\nabla_{\textstyle \frac{\partial}{\partial t}}
dF
\bigr)
(e_i)
\ = \ 
\bigl(
\nabla_{\textstyle \frac{\partial}{\partial s}}
\nabla_{e_i}
dF
\bigr)
\left({\textstyle \frac{\partial}{\partial t}}\right)
\label{eq:Ricci_f_I_1} \\ 
& = &
\bigl(
\nabla_{e_i}
\nabla_{\textstyle \frac{\partial}{\partial s}}
dF
\bigr)
\left({\textstyle \frac{\partial}{\partial t}}\right)
\ - \ 
{}^N\!R\left(
dF\bigl(e_i\bigr),\,dF({\textstyle \frac{\partial}{\partial s}})
\right) 
dF\bigl({\textstyle \frac{\partial}{\partial t}} \bigr)
\nonumber \\ 
&  & 
\quad 
(\ \because\ \,
\nabla_{{\textstyle \frac{\partial}{\partial s}}}
{\textstyle \frac{\partial}{\partial t}} 
\ = \
\nabla_{e_i}
{\textstyle \frac{\partial}{\partial t}} 
\ = \ 0 
\ )
\nonumber \\
& = &
\nabla_{e_i}
\mbox{\text Hess}_{F}
\bigl(
{\textstyle \frac{\partial}{\partial s}}
,\,
{\textstyle \frac{\partial}{\partial t}}
\bigr)
\ - \ 
{}^N\!R\left(
dF(e_i),\,dF\bigl({\textstyle \frac{\partial}{\partial s}} \bigr)
\right) 
dF\bigl({\textstyle \frac{\partial}{\partial t}} \bigr)\,. 
\nonumber 
\end{eqnarray}
Then by $(\ref{eq:I_1_1})$ and 
$(\ref{eq:Ricci_f_I_1})$, we have 
\begin{eqnarray}
\mathrm{I}_1 
& = &
\sum_{i,\,j} 
h\bigl(\nabla_{e_i}\mbox{\text Hess}_{F}
\bigl(
{\textstyle \frac{\partial}{\partial s}}
,\,
{\textstyle {\textstyle \frac{\partial}{\partial t}}}
\bigr),\,\sigma_F(e_j)\bigr)
\label{eq:I_1}\\
&  & 
\qquad 
\ - \ 
\sum_{i,\,j} 
h\bigl({}^N\!R\left(
dF(e_i),\,X
\right) 
Y,\,dF(e_j)\bigr)\,
T_F(e_i,\,e_j)
\nonumber \\ 
&  & 
\qquad
\ + \ 
\sum_{i,\,j} 
h\bigl(
\nabla_{e_i}\bigl(dF({\textstyle \frac{\partial}{\partial s}})\bigr),\,
\nabla_{e_j}\bigl(dF({\textstyle \frac{\partial}{\partial t}})\bigr)\,
\bigr)\,
T_F(e_i,\,e_j)
\nonumber 
\end{eqnarray}
In the last equality, 
we used Lemma 2 for $Z$ $=$ 
$\nabla_{e_i}\mbox{\text Hess}_{F}
\bigl(
{\textstyle \frac{\partial}{\partial s}}
,\,
{\textstyle \frac{\partial}{\partial t}}
\bigr)$ 
and $W$ $=$ $e_i$. 
On the other hand since 
\begin{eqnarray*}
\sum_{i,\,j}
g(e_i,\,e_j) 
\frac{\partial T_F(e_i,\,e_j)}{\partial t}
& = & 
\frac{\partial}{\partial t}
\left(
\sum_{i,\,j}
g(e_i,\,e_j) 
T_F(e_i,\,e_j)
\right)
\ = \ 0 
\end{eqnarray*} 
by Lemma 1 (d) and 
\begin{eqnarray*}
\frac{\partial \|dF\|^2}{\partial t}
\ = \ 
\frac{\partial}{\partial t}
\sum_{j}
h\bigl(dF(e_j),\,dF(e_j)\bigr)
\ = \ 
\sum_{j}
\frac{\partial}{\partial t}
h\bigl(dF(e_j),\,dF(e_j)\bigr)\,,
\end{eqnarray*}
we have  
\begin{eqnarray}
&  & 
 \label{eq:I_2} \\ 
\mathrm{I}_2
& = & 
\frac{1}{2}\sum_{i,\,j}
\frac{\partial}{\partial s}
\left\{
h\bigl(dF(e_i),\,dF(e_j)\bigr)
\ - \ 
\frac{1}{m} 
\|dF\|^2\,g(e_i,\,e_j)
\right\}\,
\frac{\partial T_F(e_i,\,e_j)}{\partial t}
\nonumber \\ 
& = & 
\frac{1}{2}\sum_{i,\,j}
\left\{
\frac{\partial}{\partial s}
h\bigl(dF(e_i),\,dF(e_j)\bigr)
\right\}\,
\frac{\partial T_F(e_i,\,e_j)}{\partial t}
\nonumber \\ 
&  & 
\qquad\qquad
\ - \ 
\frac{1}{2m} 
\frac{\partial \|dF\|^2}{\partial s}
\sum_{i,\,j}
g(e_i,\,e_j) 
\frac{\partial T_F(e_i,\,e_j)}{\partial t}
\nonumber \\ 
& = & 
\frac{1}{2}\sum_{i,\,j}
\left\{
\frac{\partial}{\partial s}
h\bigl(dF(e_i),\,dF(e_j)\bigr)
\right\}\,
\frac{\partial T_F(e_i,\,e_j)}{\partial t}
\nonumber \\ 
& = & 
\frac{1}{2}\sum_{i,\,j}
\left\{
\frac{\partial}{\partial s}
h\bigl(dF(e_i),\,dF(e_j)\bigr)
\right\}\,
\frac{\partial}{\partial t}
\left\{
h\bigl(dF(e_i),\,dF(e_j)\bigr)
\ - \ 
\frac{1}{m} 
\|dF\|^2\,g(e_i,\,e_j)
\right\}\,
\nonumber \\ 
& = & 
\frac{1}{2}\sum_{i,\,j}
\left\{
\frac{\partial}{\partial s}
h\bigl(dF(e_i),\,dF(e_j)\bigr)
\right\}\,
\left\{
\frac{\partial}{\partial t}
h\bigl(dF(e_i),\,dF(e_j)\bigr)
\ - \ 
\frac{1}{m} 
\frac{\partial \|dF\|^2}{\partial t}
g(e_i,\,e_j)
\right\}
\nonumber \\ 
& = & 
\frac{1}{2}\sum_{i,\,j}
\left\{
\frac{\partial}{\partial s}
h\bigl(dF(e_i),\,dF(e_j)\bigr)
\right\}\,
\left\{
\frac{\partial}{\partial t}
h\bigl(dF(e_i),\,dF(e_j)\bigr)
\right\}\,
\nonumber \\ 
&  & 
\qquad\qquad
\ - \ 
\frac{1}{2m} 
\sum_{i,\,j}
\left\{
\frac{\partial}{\partial s}
h\bigl(dF(e_i),\,dF(e_j)\bigr)
\right\}\,
g(e_i,\,e_j)\,
\frac{\partial \|dF\|^2}{\partial t}
\nonumber \\ 
& = & 
\frac{1}{2}\sum_{i,\,j}
\left\{
\frac{\partial}{\partial s}
h\bigl(dF(e_i),\,dF(e_j)\bigr)
\right\}\,
\left\{
\frac{\partial}{\partial t}
h\bigl(dF(e_i),\,dF(e_j)\bigr)
\right\}\,
\nonumber \\ 
&  & 
\qquad\qquad
\ - \ 
\frac{1}{2m} 
\sum_{i}
\left\{
\frac{\partial}{\partial s}
h\bigl(dF(e_i),\,dF(e_i)\bigr)
\right\}\,
\sum_{j}
\left\{
\frac{\partial}{\partial t}
h\bigl(dF(e_j),\,dF(e_j)\bigr)
\right\}\,
\nonumber \\ 
& = \ : & 
{\mathrm I}_3\ + \ {\mathrm I}_4 
\nonumber 
\end{eqnarray}
We have 
\begin{eqnarray}
\qquad {\mathrm I}_3 
& = & 
\frac{1}{2}\sum_{i,\,j}
\left\{
h\bigl(\nabla_{{\textstyle \frac{\partial}{\partial s}}}\bigl(dF(e_i)\bigr),\,dF(e_j)\bigr)
\ + \ 
h\bigl(dF(e_i),\,\nabla_{{\textstyle \frac{\partial}{\partial s}}}\bigl(dF(e_j)\bigr)\bigr)
\right\}\,
\label{eq:I_3} \\ 
&  & 
\qquad\qquad 
\times 
\left\{
h\bigl(\nabla_{{\textstyle \frac{\partial}{\partial t}}}\bigl(dF(e_i)\bigr),\,dF(e_j)\bigr)
\ + \ 
h\bigl(dF(e_i),\,\nabla_{{\textstyle \frac{\partial}{\partial t}}}\bigl(dF(e_j)\bigr)\bigr)
\right\}\,
\nonumber \\ 
& = & 
\sum_{i,\,j}
h\bigl(
\nabla_{{\textstyle \frac{\partial}{\partial s}}}\bigl(dF(e_i)\bigr)
,\,
dF(e_j)
\bigr)
h\bigl(
\nabla_{{\textstyle \frac{\partial}{\partial t}}}\bigl(dF(e_i)\bigr)
,\,
dF(e_j)
\bigr)
\nonumber \\ 
&  & 
\qquad\qquad
\ + \ 
\sum_{i,\,j}
h\bigl(
\nabla_{{\textstyle \frac{\partial}{\partial s}}}\bigl(dF(e_i)\bigr)
,\,
dF(e_j)
\bigr)
h\bigl(
dF(e_i)
,\,
\nabla_{{\textstyle \frac{\partial}{\partial t}}}\bigl(dF(e_j)\bigr)
\bigr)
\nonumber \\ 
& = & 
\sum_{i,\,j}
h\bigl(
\nabla_{e_i}\bigl(dF({\textstyle \frac{\partial}{\partial s}})\bigr)
,\,
dF(e_j)
\bigr)
h\bigl(
\nabla_{e_i}\bigl(dF({\textstyle \frac{\partial}{\partial t}})\bigr)
,\,
dF(e_j)
\bigr)
\nonumber \\ 
&  & 
\qquad\qquad
\ + \ 
4\,
\sum_{i,\,j}
h\bigl(
\nabla_{e_i}\bigl(dF({\textstyle \frac{\partial}{\partial s}})\bigr)
,\,
dF(e_j)
\bigr)
h\bigl(
dF(e_i)
,\,
\nabla_{e_j}\bigl(dF({\textstyle \frac{\partial}{\partial t}})\bigr)
\bigr)\,. 
\nonumber 
\end{eqnarray}
We see 
\begin{eqnarray}
\qquad {\mathrm I}_4
& = & 
\ - \ 
\frac{2}{m}\,
\sum_{i}
h\bigl(
\nabla_{{\textstyle \frac{\partial}{\partial s}}}\bigl(dF(e_i)\bigr)
,\,
dF(e_i)
\bigr)
\sum_{j}
h\bigl(
\nabla_{{\textstyle \frac{\partial}{\partial t}}}\bigl(dF(e_j)\bigr)
,\,
dF(e_j)
\bigr)
\label{eq:I_4} \\ 
& = & 
\ - \ 
\frac{2}{m}\,
\sum_{i}
h\bigl(
\nabla_{e_i}\bigl(dF({\textstyle \frac{\partial}{\partial s}})\bigr)
,\,
dF(e_i)
\bigr)
\sum_{j}
h\bigl(
\nabla_{e_j}\bigl(dF({\textstyle \frac{\partial}{\partial t}})\bigr)
,\,
dF(e_j)
\bigr)\,.
\nonumber 
\end{eqnarray}
Then by $(\ref{eq:I_1})$, $(\ref{eq:I_2})$, 
$(\ref{eq:I_3})$ and $(\ref{eq:I_4})$, we have 
\begin{eqnarray*}
\lefteqn{
\frac{1}{4}
\left.\frac{\partial^2\Phi(f_{s,\,t})}{\partial s \partial t}\right|_{s,\,t = 0}
\ = \ 
\int_M  
h\bigl(
\mbox{\text Hess}_{F}
\bigl(
{\textstyle \frac{\partial}{\partial s}}
,\,
{\textstyle \frac{\partial}{\partial t}}
\bigr),\,
\mathrm{div}_g\sigma_{f}
\bigr)
dv_g
} 
\nonumber \\ 
& &
\ + \ 
\int_M
\sum_{i,\,j} 
h\bigl(
\nabla_{e_i} \bigl(dF({\textstyle \frac{\partial}{\partial s}})\bigr),\,
\nabla_{e_j} \bigl(dF({\textstyle \frac{\partial}{\partial t}})\bigr) \bigr)\,
T_f(e_i,\,e_j)\,
dv_g
\nonumber \\ 
&  & 
\ + \ 
\int_M
\sum_{i,\,j} 
h\bigl(\nabla_{e_i} \bigl(dF({\textstyle \frac{\partial}{\partial s}})\bigr),\,df(e_j) \bigr)\,
h\bigl(\nabla_{e_i} \bigl(dF({\textstyle \frac{\partial}{\partial t}})\bigr),\,df(e_j)\bigr) 
dv_g
\nonumber \\ 
&  & 
\ + \ 
\int_M
\sum_{i,\,j} 
h\left(\nabla_{e_i} \bigl(dF({\textstyle \frac{\partial}{\partial s}})\bigr),\,df(e_j)\right)\, 
h\left(df(e_i),\,\nabla_{e_j} \bigl(dF({\textstyle \frac{\partial}{\partial t}})\bigr)\right) 
dv_g
\nonumber \\ 
&  & 
\ - \ 
\frac{2}{m}
\int_M
\sum_{i} 
h\left(\nabla_{e_i} \bigl(dF({\textstyle \frac{\partial}{\partial s}})\bigr),\,df(e_i)\right)\, 
\sum_{j} 
h\left(\nabla_{e_j} \bigl(dF({\textstyle \frac{\partial}{\partial t}})\bigr),\,df(e_j)\right)\, 
dv_g
\nonumber \\ 
&  & 
\ + \ 
\int_M
\sum_{i,\,j} 
h\bigl({}^N\!
R\left(
dF({\textstyle \frac{\partial}{\partial s}}),\,
df(e_i)\,
\right)dF({\textstyle \frac{\partial}{\partial t}}),\,
df(e_j)
\bigr)\,
T_f(e_i,\,e_j)\,
dv_g\,.
\nonumber 
\end{eqnarray*}
Integrate it over $M$ 
and let $t = 0$. 
Then using the integration by parts, 
we obtain the second variation formula.  
$\square$

\end{document}